\documentclass[psamsfonts,reqno]{amsart}
\usepackage{amssymb,eucal,graphics,latexsym}

\newcommand{\bD}{\mathbf{D}}
\newcommand{\bM}{\mathbf{M}}
\newcommand{\bMo}{\mathbf{M}_{\boldsymbol{\omega}}}
\newcommand{\SUB}{\mathbf{SUB}}
\newcommand{\CC}{\SUB(\mathcal{LO})}

\newcommand{\dnw}{\mathbin{\downarrow}}

\newcommand{\perm}{\mathfrak{S}}
\newcommand{\id}{\mathrm{id}}

\newcommand{\pI}[1]{\bigl(#1\bigr)}
\newcommand{\pII}[1]{\Bigl(#1\Bigr)}
\newcommand{\pIII}[1]{\biggl(#1\biggr)}

\newcommand{\PP}{\mathcal{P}}

\newcommand{\jh}{join-ho\-mo\-mor\-phism}
\newcommand{\mh}{meet-ho\-mo\-mor\-phism}

\newcommand{\IIdistr}{$2$-dis\-trib\-u\-tive}
\newcommand{\IIdistry}{$2$-dis\-trib\-u\-tiv\-i\-ty}
\newcommand{\jirr}{join-ir\-re\-duc\-i\-ble}

\newcommand{\jirry}{join-ir\-re\-duc\-i\-bil\-i\-ty}

\newcommand{\jsd}{join-sem\-i\-dis\-trib\-u\-tive}
\newcommand{\jsdy}{join-sem\-i\-dis\-trib\-u\-tiv\-i\-ty}

\newcommand{\mjc}{minimal nontrivial join-cover}
\newcommand{\contr}{a contradiction}

\newcommand{\CB}{\mathcal{B}}

\newcommand{\CK}{\mathcal{K}}
\newcommand{\CS}{\mathcal{S}}
\newcommand{\CX}{\mathcal{X}}
\newcommand{\VV}{\mathbf{V}}
\newcommand{\QQ}{\mathbf{Q}}
\newcommand{\DD}{\mathbin{D}}
\newcommand{\rd}[1]{[{#1}]^{\DD}}
\newcommand{\pup}[1]{\textup{(}{#1}\textup{)}}

\newcommand{\ol}[1]{\overline{#1}}

\newcommand{\fsi}[1]{\{1,\dots,#1\}}
\newcommand{\fso}[1]{\{0,\dots,#1\}}

\newcommand{\es}{\varnothing}

\newcommand{\tr}{\vartriangleleft}

\newcommand{\utr}{\trianglelefteq}
\newcommand{\gtr}{\trianglerighteq}
\newcommand{\nutr}{\ntrianglelefteq}
\newcommand{\tra}{\tr_a}

\newcommand{\utra}{\utr_a}
\newcommand{\nutra}{\nutr_a}
\newcommand{\dtr}{\mathbin{\vartriangleleft\kern-10pt {\lower
3pt\hbox{$\scriptscriptstyle\neq$}}\kern3pt}}

\newcommand{\bel}{\utr_{\mathrm{w}}}

\newcommand{\gf}{\varphi}
\newcommand{\gfa}{\varphi_a}
\newcommand{\gy}{\psi}
\newcommand{\gya}{\psi_a}

\newcommand{\set}[1]{\{{#1}\}}
\newcommand{\setm}[2]{\set{{#1}\mid{#2}}}
\newcommand{\seq}[1]{\langle#1\rangle}
\newcommand{\famm}[2]{\seq{{#1}\nobreak\mid\nobreak{#2}}}

\newcommand{\St}{\textup{(S)}}
\newcommand{\Ud}{\textup{(U)}}
\newcommand{\Bo}{\textup{(B)}}
\newcommand{\El}{\textup{(E)}}
\newcommand{\Pe}{\textup{(P)}}
\newcommand{\HS}{\textup{(HS)}}

\newcommand{\Ht}[1]{\textup{(H$_{#1}$)}}

\newcommand{\Stj}{\textup{(S$_{\mathrm{j}}$)}}
\newcommand{\Udj}{\textup{(U$_{\mathrm{j}}$)}}
\newcommand{\Boj}{\textup{(B$_{\mathrm{j}}$)}}
\newcommand{\Elj}{\textup{(E$^{\Sigma}$)}}
\newcommand{\Pej}{\textup{(P$^{\Sigma}$)}}
\newcommand{\HSj}{\textup{(HS$^{\Sigma}$)}}

\newcommand{\two}{\mathbf{2}}

\newcommand{\hL}{{\widehat{L}}}

\newcommand{\into}{\hookrightarrow}
\newcommand{\onto}{\twoheadrightarrow}

\DeclareMathOperator{\J}{J}

\newcommand{\Var}{\mathbf{V}}
\newcommand{\QVar}{\mathbf{Q}}

\newcommand{\Co}{\mathbf{Co}}

\newcommand{\Csub}{\mathbf{Csub}}

\numberwithin{equation}{section}

\theoremstyle{plain}
\newtheorem{lemma}{Lemma}[section]
\newtheorem{theorem}[lemma]{Theorem}
\newtheorem{proposition}[lemma]{Proposition}
\newtheorem{corollary}[lemma]{Corollary}

\newtheorem*{stat}{\name}
\newcommand{\name}{testing}

\theoremstyle{definition}
\newtheorem{definition}[lemma]{Definition}
\newtheorem{notation}[lemma]{Notation}
\newtheorem{example}[lemma]{Example}
\newtheorem{problem}{Problem}

\theoremstyle{remark}
\newtheorem{remark}[lemma]{Remark}

\newtheorem*{noteadd}{Note added}

\newenvironment{all}[1]{\renewcommand{\name}{#1}\begin{stat}}
                          {\end{stat}}

\newcommand{\qedc}{{\qed}~{\rm Claim~{\theclaim}.}}
\newcommand{\qedsc}{{\qed}~{\rm Claim.}}

\begin{document}

\title[Order-convex subsets of chains]%
{Sublattices of lattices of order-convex sets, III.\\
The case of totally ordered sets}

\author[M.~Semenova]{Marina Semenova}
\address[M.~Semenova]{Institute of Mathematics of
the Siberian Branch of RAS\\
Acad. Koptyug prosp. 4\\
630090 Novosibirsk\\
Russia}
\email{semenova@math.nsc.ru}

\author[F.~Wehrung]{Friedrich Wehrung}
\address[F.~Wehrung]{CNRS, UMR 6139\\
D\'epartement de Math\'ematiques\\
Universit\'e de Caen\\
14032 Caen Cedex\\
France}
\email{wehrung@math.unicaen.fr}
\urladdr{http://www.math.unicaen.fr/\~{}wehrung}
\date{\today}

\subjclass[2000]{Primary: 06B05, 06B20, 06B15, 06A05, 08C15.
Secondary: 05B25}
\keywords{Lattice, embedding, poset, chain, order-convex, variety,
\jirr, join-seed}

\thanks{The first author was partially supported by INTAS grant no.
YSF: 2001/1-65. The authors were partially supported by GA CR grant
no. 201/00/0766 and by institutional grant MSM:J13/98:1132000007a}

\begin{abstract}
For a partially ordered set $P$, let $\Co(P)$ denote the lattice of
all order-convex subsets of $P$.
For a positive integer $n$, we denote by $\CC$ (resp.,
$\SUB(n)$) the class of all lattices that can be embedded into a
lattice of the form
   \[
   \prod_{i\in I}\Co(T_i),
   \]
where $\famm{T_i}{i\in I}$ is a family of \emph{chains} (resp.,
chains with at most $n$ elements). We prove the following results:
\begin{itemize}
\item[(1)] Both classes $\CC$ and $\SUB(n)$, for any positive integer
$n$, are locally finite, finitely based varieties of lattices, and we
find finite equational bases of these varieties.

\item[(2)] The variety $\CC$ is the quasivariety join of
all the varieties $\SUB(n)$, for $1\leq n<\omega$, and it has only
countably many subvarieties. We classify these varieties, together with
all the finite subdirectly irreducible members of $\CC$.

\item[(3)] Every finite subdirectly irreducible member of $\CC$
is projective within $\CC$, and every subquasivariety of
$\CC$ is a variety.

\end{itemize}
\end{abstract}

\maketitle

\section{Introduction}\label{S:Intro}

For a partially ordered set (from now on \emph{poset}) $(P,\utr)$, a
subset $X$ of $P$ is \emph{order-convex}, if $x\utr z\utr y$ and
$\set{x,y}\subseteq X$ implies that $z\in X$, for all $x$, $y$,
$z\in P$. The lattices of the form $\Co(P)$ have been characterized
by G. Birkhoff and M.\,K. Bennett in \cite{BB}. In M. Semenova and
F. Wehrung \cite{SeWe1}, the authors solve a problem stated in K.\,V.
Adaricheva, V.\,A. Gorbunov, and V.\,I. Tumanov \cite{AGT}, by proving
the following result.

\begin{all}{Theorem~1}
The class $\SUB$ of all lattices that can be embedded into some
lattice of the form $\Co(P)$ forms a variety, defined by three
identities, \St, \Ud, and \Bo.
\end{all}

In M. Semenova and F. Wehrung \cite{SeWe2}, this result is extended
to special classes of posets $P$:

\begin{all}{Theorem~2}
For a positive integer $n$, the class $\SUB_n$ of all lattices that
can be embedded into some lattice of the form $\Co(P)$, where $P$
is a poset of length at most~$n$, is a variety, defined by the
identities \St, \Ud, \Bo, together with new identities \Ht{n} and
\Ht{k,n+1-k}, for $1\leq k\leq n$.
\end{all}

In the present paper, we extend these results to sublattices of
products of lattices of convex subsets of \emph{chains} (i.e.,
totally ordered sets), thus solving a problem of~\cite{SeWe1}. More
specifically, we denote by $\CC$ (resp., $\SUB(n)$) the class of all
lattices that can be embedded into a lattice of the form
   \[
   \prod_{i\in I}\Co(T_i),
   \]
where $\famm{T_i}{i\in I}$ is a family of chains (resp.,
chains with at most $n$ elements). We prove the following results:
\begin{itemize}
\item[(1)] Both classes $\CC$ and $\SUB(n)$ are finitely based
varieties of lattices, for any positive integer $n$. Moreover,
$\SUB(n+1)=\CC\cap\SUB_n$ (Theorems~\ref{T:EmbThm} and
\ref{T:SUB(n)}).

\item[(2)] By using a result of V. Slav\'\i k \cite{Slav}, we prove
that the variety $\CC$ is locally finite (Theorem~\ref{T:CLocFin}).

\item[(3)] The variety $\CC$ is the quasivariety join of
all the varieties $\SUB(n)$, for $1\leq n<\omega$
(Corollary~\ref{C:QVjoin}), and every proper
subvariety of $\CC$ is finitely generated
(Corollary~\ref{C:ClassSI}).

\item[(4)] The only proper subvarieties of $\CC$ are those between
$\SUB(n)$ and $\SUB(n+1)$ for some natural number $n$
(Theorem~\ref{T:ClassSI}).

\item[(5)] We classify all finite subdirectly irreducible members
of $\CC$, and we describe exactly the lattice of all subvarieties
of $\CC$ (Theorem~\ref{T:ClassSI} to Corollary~\ref{C:subSUB(n)}).

\item[(6)] All finite subdirectly irreducible members of
$\CC$ are projective within $\CC$ (Theorem~\ref{T:SIProj}), and every
subquasivariety of $\CC$ is a variety (Theorem~\ref{T:QVar2Var}).
\end{itemize}

The main technical result towards the proof that $\CC$ is a variety
is that the reflexive closure of the join-dependency relation $\DD$
is \emph{transitive}, in any member of $\CC$ with `enough' \jirr\
elements (Corollary~\ref{C:Trans}). This may be
viewed as an analogue, for certain \jsd\ lattices, of the transitivity
of perspectivity proved by von~Neumann in continuous geometries,
see \cite{Neum}.

We refer the reader to our papers \cite{SeWe1,SeWe2} for unexplained
notation and terminology. In particular,
the identities \St, \Ud, and \Bo, together with their \jirr\
translations \Stj, \Udj, and \Boj, and tools such as Stirlitz
tracks or the Udav-Bond partition, are defined in
\cite{SeWe1}. The identities \Ht{n} and \Ht{m,n}, their \jirr\
translations, and bi-Stirlitz tracks are defined in \cite{SeWe2}.
We shall often use the trivial fact that $\Co(P,\utr)=\Co(P,\gtr)$, for
any poset $(P,\utr)$, where $\gtr$ denotes the converse order of $\utr$.

The join-dependency relation on a lattice $L$, see R. Freese, J.
Je\v{z}ek, and J.\,B. Nation \cite{FJN}, is defined on the set
$\J(L)$ of all \jirr\ elements of $L$, and it is written
$\DD_L$, or $\DD$ if $L$ is understood from the context. For
$a\in\J(L)$, we write, as in \cite{SeWe1,SeWe2},
   \[
   \rd{a}=\setm{x\in\J(L)}{a\DD x}.
   \]

\section{Join-seeds and more minimal covers}\label{S:JoinSeeds}

We recall from \cite{SeWe2} the following definition:

\begin{definition}\label{D:GoodGen}
A subset $\Sigma$ of a lattice $L$ is a \emph{join-seed}, if the
following statements hold:
\begin{enumerate}
\item $\Sigma\subseteq\J(L)$;

\item every element of $L$ is a join of elements of $\Sigma$;

\item for all $p\in\Sigma$ and all $a$, $b\in L$ such that
$p\leq a\vee b$ and $p\nleq a,b$, there are $x\leq a$ and $y\leq b$
both in $\Sigma$ such that $p\leq x\vee y$ is minimal in $x$ and $y$.
\end{enumerate}
\end{definition}

Two important examples of join-seeds are provided by the following
lemma, see~\cite{SeWe2}.

\begin{lemma}\label{L:JoinSeeds}
Any of the following assumptions implies that the subset $\Sigma$ is
a join-seed of the lattice $L$:
\begin{enumerate}
\item $L=\Co(P)$ and $\Sigma=\setm{\set{p}}{p\in P}$, for some
poset $P$.

\item $L$ is a dually \IIdistr, complete, lower continuous,
finitely spatial lattice, and $\Sigma=\J(L)$.
\end{enumerate}
\end{lemma}

\begin{lemma}\label{L:MinUB}
Let $L$ be a lattice satisfying \Bo, let $\Sigma$ be a join-seed of
$L$, let $p\in\Sigma$, let $x$, $y\in\rd{p}$. If the inequality
$p\leq x\vee y$ holds, then it is minimal in both $x$ and~$y$.
\end{lemma}

\begin{proof}
{}From the assumption that $x$, $y\in\rd{p}$, it follows
that $p\nleq x,y$. Since $p\leq x\vee y$ and $\Sigma$ is a join-seed of
$L$, there are $u\leq x$ and $v\leq y$ in $\Sigma$ such that the
inequality $p\leq u\vee v$ holds and is minimal in both $u$ and $v$.
Furthermore, by the definition of the $\DD$ relation and since
$\Sigma$ is a join-seed of $L$, there are $x'$, $y'\in\Sigma$ such
that both inequalities $p\leq x\vee x'$ and $p\leq y\vee y'$ hold
and are minimal in $x$, $x'$, $y$, $y'$. By applying \Boj\ to the
inequalities $p\leq x\vee x',u\vee v$ and by observing that
$p\nleq x,v$, we obtain that $p\leq x'\vee u$. Since $u\leq x$ and
the inequality $p\leq x\vee x'$ is minimal in~$x$, we obtain that
$u=x$. Similarly, $v=y$.
\end{proof}

\begin{lemma}\label{L:rdaAnt}
Let $L$ be a lattice satisfying \Bo, let $\Sigma$ be a join-seed of
$L$. Then $\rd{p}\cap\Sigma$ is an antichain of $L$, for any
$p\in\Sigma$.
\end{lemma}

\begin{proof}
Let $x$, $y\in\rd{p}$. Since $\Sigma$ is a join-seed of $L$, there
are $x'$, $y'\in\Sigma$ such that both inequalities $p\leq x\vee x'$
and $p\leq y\vee y'$ are \mjc s. Observe that $p\nleq x,x',y,y'$. If
$x\leq y$, then, since $p\nleq y=x\vee y$ and $L$ satisfies \Boj, the
inequality $p\leq x\vee y'$ holds. Since $x\leq y$ and the inequality
$p\leq y\vee y'$ is minimal in $y$, we obtain that $x=y$.
\end{proof}

\section{The identity \El}\label{S:Elash}

Let \El\ be the following identity in the variables $x$, $a$, $b_0$,
$b_1$, $b_2$:
   \begin{multline*}
   x\wedge\bigwedge_{i<3}(a\vee b_i)=
   \bigvee_{i<3}\left[x\wedge b_i\wedge\bigwedge_{j\neq i}
   (a\vee b_j)\right]\\
   \vee\bigvee_{\sigma\in\perm_3}\bigl[
   x\wedge(a\vee b^*_{0,\sigma})
   \wedge
   (a\vee b^*_{1,\sigma})
   \wedge(a\vee b_{\sigma(2)})
   \bigr],
   \end{multline*}
where we denote by $\perm_3$ the group of all permutations of
$\set{0,1,2}$ and we put
   \begin{align}
   b^*_{0,\sigma}&=b_{\sigma(0)}\wedge(x\vee b_{\sigma(1)}),
   \label{Eq:b*0s}\\
   b^*_{1,\sigma}&=b_{\sigma(1)}\wedge(x\vee b_{\sigma(2)})\wedge
   (b_{\sigma(0)}\vee b_{\sigma(2)})\label{Eq:b*1s},
   \end{align}
for all $\sigma\in\perm_3$.

We now introduce a lattice-theoretical axiom, the \emph{\jirr\
interpretation of \El}, that we will denote by \Elj.

\begin{definition}\label{D:Elj}
For a lattice $L$ and a subset $\Sigma$ of $\J(L)$, we say that $L$
satisfies $\Elj$, if for all elements $x$, $a$, $b_0$, $b_1$, and
$b_2$ of $\Sigma$, if the inequality $x\leq a\vee b_i$ is
a \mjc, for every $i<3$, then there exists
$\sigma\in\perm_3$ such that
$b_{\sigma(0)}\leq x\vee b_{\sigma(1)}\leq x\vee b_{\sigma(2)}$ and
$b_{\sigma(1)}\leq b_{\sigma(0)}\vee b_{\sigma(2)}$.
\end{definition}

The geometrical meaning of \Elj\ is illustrated on Figure~1. The lines
of that figure represent the ordering of the either the poset $P$ or its dual
(and not the ordering of $L$) in case $L=\Co(P,\utr)$.
For example, the left half of Figure~1 represents (up to dualization
of $\utr$) the relations $a\utr x\utr b_i$, for $i<3$, so that the
inequality $\set{x}\leq\set{a}\vee\set{b_i}$ holds in~$L$. Similar
conventions hold for Figures~2 and~3.

\begin{figure}[htb]
\includegraphics{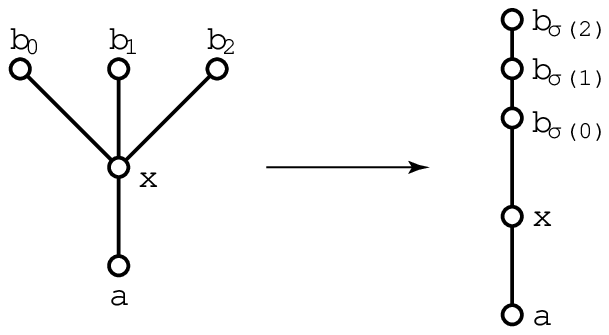}
\caption{Illustrating \Elj}
\end{figure}

\begin{lemma}\label{L:El2Elj}
Let $L$ be a lattice, let $\Sigma$ be a subset of $\J(L)$. Then the
following statements hold:
\begin{enumerate}
\item If $L$ satisfies \El, then $L$ satisfies \Elj.

\item If $\Sigma$ is a join-seed of $L$ and $L$ satisfies both \Bo\
and \Elj, then $L$ satisfies~\El.
\end{enumerate}
\end{lemma}

\begin{proof}
(i)
Suppose that $x$, $a$, $b_0$, $b_1$,
$b_2\in\Sigma$ satisfy the premise of \Elj. Since~$x$ is \jirr\ and
$x\nleq b_i$, for all $i<3$, we obtain, by applying the
identity~\El\ and using the notation introduced in \eqref{Eq:b*0s}
and \eqref{Eq:b*1s}, that there exists $\sigma\in\perm_3$ such that
both inequalities
$x\leq a\vee b^*_{0,\sigma},a\vee b^*_{1,\sigma}$ hold. Since
$b^*_{i,\sigma}\leq b_{\sigma(i)}$, it follows from
the minimality of $b_{\sigma(i)}$ in the inequality
$x\leq a\vee b_{\sigma(i)}$ that $b^*_{i,\sigma}=b_{\sigma(i)}$, for
all $i<2$. Therefore,
$b_{\sigma(0)}\leq x\vee b_{\sigma(1)}\leq x\vee b_{\sigma(2)}$ and
$b_{\sigma(1)}\leq b_{\sigma(0)}\vee b_{\sigma(2)}$.

(ii) Let $c$ (resp., $d$) denote the left hand side
(resp., right hand side) of the identity \El. Since $d\leq c$ holds
in any lattice, it suffices to prove that $c\leq d$. Let
$p\in\Sigma$ with $p\leq c$, we prove that $p\leq d$. If $p\leq a$,
then $p\leq x\wedge a\leq d$. If $p\leq b_i$, for some $i<3$, then
$p\leq x\wedge b_i\wedge\bigwedge_{j\neq i}(a\vee b_j)\leq d$.

Suppose from now on that $p\nleq a$ and $p\nleq b_i$, for all $i<3$.
Since $p\leq a\vee b_i$ and $\Sigma$ is a join-seed of $L$, there
are $u_i\leq a$ and $v_i\leq b_i$ in $\Sigma$ such that the
inequality $p\leq u_i\vee v_i$ is a \mjc, for all $i<3$. In particular,
$u_i$, $v_i\in\rd{p}$. Put
$u=u_0$, and let $i<3$. By applying \Boj\ to the inequalities
$p\leq u\vee v_0,u_i\vee v_i$ and observing that $p\nleq a$ (thus
$p\nleq u\vee u_i$), we obtain the inequality $p\leq u\vee v_i$.
Furthermore, by Lemma~\ref{L:MinUB}, this inequality is minimal in
both $u$ and $v_i$. Hence, by \Elj, there exists $\sigma\in\perm_3$
such that
$v_{\sigma(0)}\leq p\vee v_{\sigma(1)}\leq p\vee v_{\sigma(2)}$
and $v_{\sigma(1)}\leq v_{\sigma(0)}\vee v_{\sigma(2)}$. Therefore,
by putting
   \begin{align*}
   v^*_{0,\sigma}&=v_{\sigma(0)}\wedge(p\vee v_{\sigma(1)}),\\
   v^*_{1,\sigma}&=v_{\sigma(1)}\wedge(p\vee v_{\sigma(2)})\wedge
   (v_{\sigma(0)}\vee v_{\sigma(2)}),
   \end{align*}
we obtain the equalities $v^*_{0,\sigma}=v_{\sigma(0)}$ and
$v^*_{1,\sigma}=v_{\sigma(1)}$, and the inequalities
   \[
   p\leq x\wedge(u\vee v^*_{0,\sigma})\wedge(u\vee v^*_{1,\sigma})
   \wedge(u\vee v_{\sigma(2)})\leq d.
   \]
Since every element of $L$ is a join of elements of $\Sigma$, the
inequality $c\leq d$ follows.
\end{proof}

\begin{corollary}\label{C:El2Elj}
The lattice $\Co(T)$ satisfies the
identity \El, for any chain $(T,\utr)$.
\end{corollary}

\begin{proof}
We apply Lemma~\ref{L:El2Elj} to $L=\Co(T)$ together with the
join-seed $\Sigma=\setm{\set{p}}{p\in T}$. Let $x$, $a$, $b_0$,
$b_1$, $b_2\in T$ such that the inequality
$\set{x}\leq\set{a}\vee\set{b_i}$ is a \mjc,
for all $i<3$. Since $\Co(T,\utr)=\Co(T,\gtr)$, we may assume without
loss of generality that $a\tr x\tr b_0$, thus $x\tr b_i$, for
all $i<3$. Since~$T$ is a chain, there exists $\sigma\in\perm_3$
such that $b_{\sigma(0)}\utr b_{\sigma(1)}\utr b_{\sigma(2)}$, whence
   \[
   \set{b_{\sigma(0)}}\leq\set{x}\vee\set{b_{\sigma(1)}}\leq
   \set{x}\vee\set{b_{\sigma(2)}}\text{ and }
   \set{b_{\sigma(1)}}\leq\set{b_{\sigma(0)}}\vee\set{b_{\sigma(2)}}.
   \]
Hence $\Co(T)$ satisfies \Elj.
Since $\Co(T)$ satisfies \Bo\ (see \cite{SeWe1}) and $\Sigma$ is a
join-seed of $\Co(T)$, it follows from Lemma~\ref{L:El2Elj} that
$\Co(T)$ satisfies \El.
\end{proof}

\begin{lemma}\label{L:IVP}
Let $L$ be a \jsd\ lattice satisfying the identity \El, let
$a$, $x\in\J(L)$ and $b_0$, $b_1$, $b_2\in\J(L)$ be distinct such
that $x\leq a\vee b_i$ is a \mjc, for all
$i<3$. Then $a\vee b_0\leq a\vee b_1\leq a\vee b_2$ implies that
$a\vee b_0<a\vee b_1<a\vee b_2$ and $b_1\leq b_0\vee b_2$.
\end{lemma}

\begin{proof}
Let $i$, $j$ be distinct in $\set{0,1,2}$. If $a\vee b_i=a\vee b_j$,
then, by the \jsdy\ of $L$, $x\leq a\vee b_i=a\vee(b_i\wedge b_j)$;
it follows from the minimality assumption on $b_i$ that $b_i\leq b_j$.
Similarly, $b_j\leq b_i$, whence $b_i=b_j$, \contr.
Thus we have obtained the inequalities
   \begin{equation}\label{Eq:ab0b1b2<}
   a\vee b_0<a\vee b_1<a\vee b_2.
   \end{equation}
On the other hand, it follows from Lemma~\ref{L:El2Elj} that
there exists $\sigma\in\perm_3$ such that the inequalities
   \begin{gather}
   x\vee b_{\sigma(0)}\leq x\vee b_{\sigma(1)}
   \leq x\vee b_{\sigma(2)},\label{Eq:xb0b1b2<}\\
   b_{\sigma(1)}\leq b_{\sigma(0)}\vee b_{\sigma(2)}\label{Eq:bs102}
   \end{gather}
hold. {}From \eqref{Eq:xb0b1b2<} it follows that
$a\vee b_{\sigma(0)}\leq a\vee b_{\sigma(1)}\leq a\vee b_{\sigma(2)}$,
thus, by \eqref{Eq:ab0b1b2<}, $\sigma$ is the identity. The
conclusion follows from \eqref{Eq:xb0b1b2<} and \eqref{Eq:bs102}.
\end{proof}

\section{The identity \Pe}\label{S:Pech}

Let \Pe\ be the following identity in the variables $a$, $b$, $c$,
$d$, $b_0$, $b_1$:
   \begin{align*}
   a\wedge(b'\vee c)\wedge(c\vee d)=&
   \bigl(a\wedge b'\wedge(c\vee d)\bigr)\vee
   \bigl(a\wedge d\wedge(b'\vee c)\bigr)\\
   &\vee\biggl[
   a\wedge\Bigl(\bigl(b'\wedge(a\vee d)\bigr)\vee c\Bigr)
   \wedge(c\vee d)\biggr]\\
   &\vee\bigvee_{i<2}\biggl[a\wedge(b_i\vee c)\wedge
   \Bigl(\bigl(b'\wedge(a\vee b_i)\wedge(b_i\vee d)\bigr)\vee c\Bigr)
   \wedge(c\vee d)\biggr],
   \end{align*}
where we put $b'=b\wedge(b_0\vee b_1)$.

We now introduce a lattice-theoretical axiom, the \emph{\jirr\
interpretation of \Pe}, that we will denote by \Pej.

\begin{definition}\label{D:Pej}
For a lattice $L$ and a subset $\Sigma$ of $\J(L)$, we say that $L$
satisfies $\Pej$, if for all elements $a$, $b$, $c$, $d$, $b_0$,
$b_1$ in $\Sigma$, if both inequalities $a\leq b\vee c,c\vee d$ are
\mjc s and $b\leq b_0\vee b_1$, then either
$b\leq a\vee d$ or there exists $i<2$ such that $a\leq b_i\vee c$
and $b\leq a\vee b_i,b_i\vee d$.
\end{definition}

The geometrical meaning of \Pej\ is illustrated on Figure~2. Horizontal
lines are meant to suggest that ``no side is chosen yet''. For example,
the non-horizontal lines in the left half of Figure~2 represent various
inequalities such as $c\utr a\utr d$ and $c\utr a\utr b$
(in case $L=\Co(P,\utr)$), while the
horizontal line represents the inequalities $b_{1-i}\utr b\utr b_i$, for
some $i<2$. A similar convention applies to Figure~3.

\begin{figure}[htb]
\includegraphics{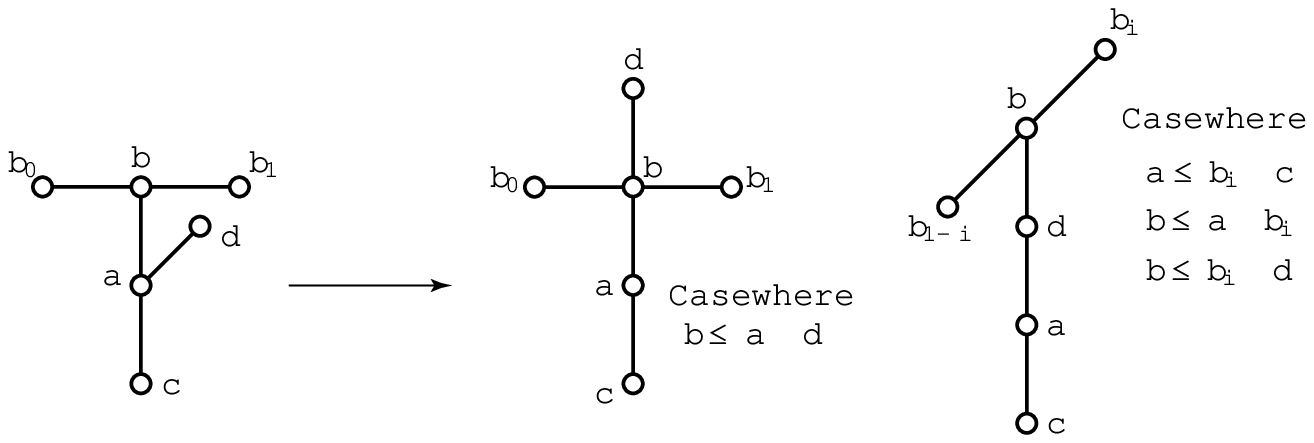}
\caption{Illustrating \Pej}
\end{figure}

\begin{lemma}\label{L:Pe2Pej}
Let $L$ be a lattice, let $\Sigma$ be a subset of $\J(L)$. Then the
following statements hold:
\begin{enumerate}
\item If $L$ satisfies \Pe, then $L$ satisfies \Pej.

\item If $\Sigma$ is a join-seed of $L$ and $L$ satisfies both \Bo\
and \Pej, then $L$ satisfies~\Pe.
\end{enumerate}
\end{lemma}

\begin{proof}
(i) Let $a$, $b$, $c$, $d$, $b_0$, $b_1\in\Sigma$
satisfy the premise of \Pej. Observe that $b\wedge(b_0\vee b_1)=b$,
thus the left hand side of the identity \Pe\ computed with these
parameters equals $a$. Since $a\nleq b,d$ and $a$ is \jirr, either
$a\leq\bigl(b\wedge(a\vee d)\bigr)\vee c$ or $a\leq b_i\vee c$ and
$a\leq\bigl(b\wedge(a\vee b_i)\wedge(b_i\vee d)\bigr)\vee c$, for
some $i<2$. In the first case, from the fact that the cover
$a\leq b\vee c$ is minimal in $b$ it follows that $b\leq a\vee d$ in the
first case, and $b\leq a\vee b_i,b_i\vee d$ in the second case.

(ii) Let $e$ (resp., $f$) denote the left hand side
(resp., right hand side) of the identity \Pe. Let $p\in\Sigma$ such
that $p\leq e$, we prove that $p\leq f$. If either $p\leq c$ or
$p\leq b'$ or $p\leq d$ this is obvious, so suppose, from now on,
that $p\nleq c,b',d$. Since $\Sigma$ is a join-seed of $L$, there
are $u\leq b'$ together with $v,v'\leq c$ and $w\leq d$ in $\Sigma$
such that both inequalities
   \begin{align}
   p&\leq u\vee v,\label{Eq:pluv}\\
   p&\leq v'\vee w\label{Eq:plv'w}
   \end{align}
are \mjc s. In particular, $u$, $v$, $v'$,
$w\in\rd{p}$. Furthermore, by applying \Boj\ to the inequalities
\eqref{Eq:pluv} and \eqref{Eq:plv'w} and observing that
$p\nleq v\vee v'$ (because $p\nleq c$), we obtain the inequality
   \begin{equation}\label{Eq:plvw}
   p\leq v\vee w.
   \end{equation}
Furthermore, it follows from Lemma~\ref{L:MinUB} that
\eqref{Eq:plvw} is a \mjc. Since $\Sigma$ is a join-seed of $L$,
there are $u_i\leq b_i$ in $\Sigma\cup\set{0}$, for $i<2$, such that
$u\leq u_0\vee u_1$. Suppose first that $u_0$, $u_1\in\Sigma$.
Since $L$ satisfies \Pej, either
   \begin{equation}\label{Eq:vluw}
   u\leq p\vee w
   \end{equation}
or
   \begin{equation}\label{Eq:notvluw}
   p\leq u_i\vee v\text{ and }u\leq p\vee u_i,u_i\vee w,
   \text{ for some }i<2.
   \end{equation}
The conclusion \eqref{Eq:notvluw} also holds if $u_j=0$, for some
$j<2$, because $u\leq u_{1-j}$.

If \eqref{Eq:vluw} holds, then
   \[
   p\leq a\wedge\Bigl(\bigl(u\wedge(p\vee w)\bigr)\vee v\Bigr)
   \wedge(v\vee w)\leq f.
   \]
If \eqref{Eq:notvluw} holds, then
   \[
   p\leq a\wedge(u_i\vee v)\wedge\Bigl(
   \bigl(u\wedge
   (p\vee u_i)\wedge(u_i\vee w)\bigr)\vee v\Bigr)
   \wedge(v\vee w)\leq f.
   \]
Since every element of $L$ is a join of elements of $\Sigma$, the
inequality $e\leq f$ follows. Since $f\leq e$ holds in any lattice,
we obtain that $e=f$.
\end{proof}

\begin{corollary}\label{C:Pe2Pej}
The lattice $\Co(T)$ satisfies \Pe, for every chain $(T,\utr)$.
\end{corollary}

\begin{proof}
We apply Lemma~\ref{L:Pe2Pej} to $L=\Co(T)$ together with the
join-seed $\Sigma=\setm{\set{p}}{p\in T}$. Let $a$, $b$, $c$, $d$,
$b_0$, $b_1\in T$ such that both inequalities
$\set{a}\leq\set{b}\vee\set{c},\set{c}\vee\set{d}$ are \mjc s and
$\set{b}\leq\set{b_0}\vee\set{b_1}$. Since $\Co(T,\utr)=\Co(T,\gtr)$,
we may assume without loss of generality that $c\tr a\tr b,d$.
Furthermore, from
$\set{b}\leq\set{b_0}\vee\set{b_1}$ it follows that there exists $i<2$
such that $b\utr b_i$. Since $T$ is a chain, either $b\utr d$ or
$d\utr b$. In the first case, $\set{b}\leq\set{a}\vee\set{d}$. In
the second case, $\set{a}\leq\set{b_i}\vee\set{c}$ and
$\set{b}\leq\set{a}\vee\set{b_i},\set{b_i}\vee\set{d}$.

Hence $\Co(T)$ satisfies \Pej. By Lemma~\ref{L:Pe2Pej}, $\Co(T)$
satisfies \Pe.
\end{proof}

\section{The identity \HS}\label{S:HypSt}

Let \HS\ be the following identity in the variables $a$, $b$, $c$,
$b_0$, $b_1$:
   \begin{align*}
   a\wedge(b'\vee c)=&(a\wedge b')\vee
   \bigvee_{i<2}\Bigl[a\wedge\bigl((b\wedge b_i)\vee c\bigr)\Bigr]\\
   &\vee\bigvee_{i<2}\biggl[
   a\wedge\Bigl(\bigl(b'\wedge(a\vee b_i)\bigr)\vee c\Bigr)
   \wedge(b_i\vee c)\wedge(b\vee b_{1-i})\biggr]\\
   &\vee\bigvee_{i<2}\biggl[
   a\wedge\Bigl(\bigl(b'\wedge(a\vee b_i)\bigr)\vee c\Bigr)
   \wedge(b_0\vee c)\wedge(b_1\vee c)\biggr],
   \end{align*}
where we put $b'=b\wedge(b_0\vee b_1)$. Since the right hand side
of \HS\ lies obviously below the right hand side of the identity
\St\ while the left hand sides are the same, we obtain immediately
the following result.

\begin{lemma}\label{L:HS2St}
The identity \HS\ implies the Stirlitz identity \St.
\end{lemma}

As observed in \cite{SeWe1}, \St\ implies both \jsdy\ and dual
\IIdistry. Therefore, we obtain the following consequence.

\begin{lemma}\label{L:HS2SDj}
The identity \HS\ implies both \jsdy\ and dual
\IIdistry.
\end{lemma}

We now introduce a lattice-theoretical axiom, the \emph{\jirr\
interpretation of \HS}, that we will denote by \HSj.

\begin{definition}\label{D:HSj}
For a lattice $L$ and a subset $\Sigma$ of $\J(L)$, we say that $L$
satisfies $\HSj$, if for all elements $a$, $b$, $c$, $b_0$,
$b_1$ in $\Sigma$, if $a\neq b$, the inequality $a\leq b\vee c$ is
minimal in $b$, and $b\leq b_0\vee b_1$ is a nontrivial join-cover,
then there exists $i<2$ such that $b\leq a\vee b_i$ and either
$a\leq b_i\vee c,b\vee b_{1-i}$ or $a\leq b_0\vee c,b_1\vee c$.
\end{definition}

The geometrical meaning of \HSj\ is illustrated on Figure~3.

\begin{figure}[htb]
\includegraphics{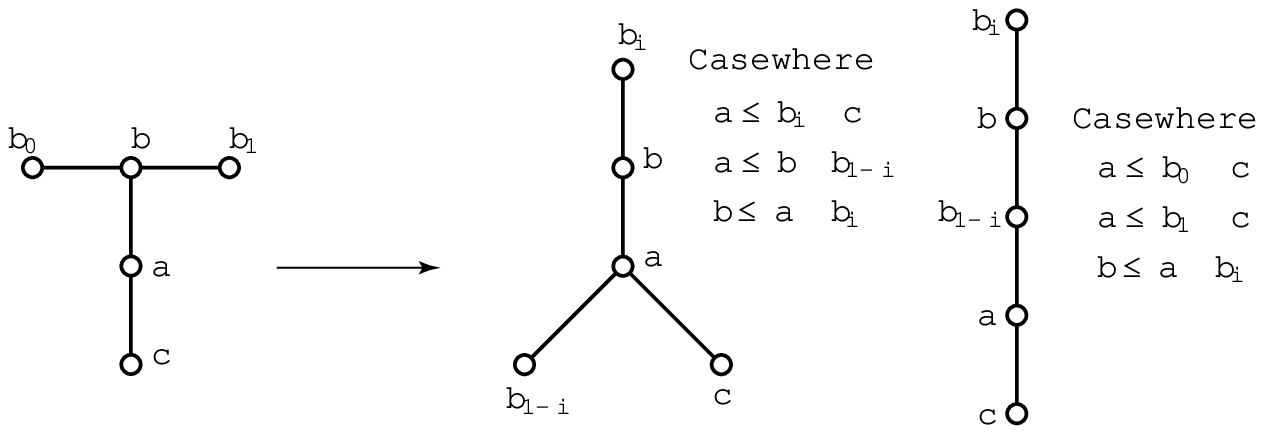}
\caption{Illustrating \HSj}
\end{figure}

\begin{lemma}\label{L:HS2HSj}
Let $L$ be a lattice, let $\Sigma$ be a subset of $\J(L)$. Then the
following statements hold:
\begin{enumerate}
\item If $L$ satisfies \HS, then $L$ satisfies \HSj.

\item If $\Sigma$ is a join-seed of $L$ and $L$ satisfies \HSj, then
$L$ satisfies~\HS.
\end{enumerate}
\end{lemma}

\begin{proof}
(i) Let $a$, $b$, $c$, $b_0$, $b_1\in\Sigma$
satisfy the premise of \HSj. Observe that $b'=b\wedge(b_0\vee b_1)=b$
and $a\wedge(b'\vee c)=a$. Since $a\leq b\vee c$ is minimal in $b$
and $b\wedge b_i<b$, it follows from the \jirry\ of $a$ that there
exists $i<2$ such that one of the following inequalities holds:
   \begin{align*}
   a&\leq\Bigl(\bigl(b\wedge(a\vee b_i)\bigr)\vee c\Bigr)
   \wedge(b_i\vee c)\wedge(b\vee b_{1-i}),\\
   a&\leq\Bigl(\bigl(b\wedge(a\vee b_i)\bigr)\vee c\Bigr)
   \wedge(b_0\vee c)\wedge(b_1\vee c).
   \end{align*}
{}From the minimality of $b$ in $a\leq b\vee c$ it follows that
$b\leq a\vee b_i$. Furthermore, in the first case
$a\leq b_i\vee c,b\vee b_{1-i}$ while in the second case
$a\leq b_0\vee c,b_1\vee c$.

(ii) Let $d$ (resp., $e$) denote the left hand side
(resp., right hand side) of the identity \HS. Let $p\in\Sigma$ such
that $p\leq d$, we prove that $p\leq e$. If $p\leq b'$ then
$p\leq d\wedge b'=a\wedge b'$, if $p\leq c$ then $p\leq a\wedge c$, in
both cases $p\leq e$. Suppose from now on that
$p\nleq b',c$. Since $\Sigma$ is a join-seed of $L$, there are
$u\leq b'$ and $v\leq c$ in $\Sigma$ such that $p\leq u\vee v$ is a
\mjc. If $u\leq b_i$, for some $i<2$, then $u\leq b\wedge b_i$,
whence
   \[
   p\leq a\wedge(u\vee v)\leq a\wedge\bigl((b\wedge b_i)\vee c\bigr)
   \leq e.
   \]
Suppose from now on that $u\nleq b_0,b_1$. Since $\Sigma$ is a
join-seed of $L$, there are $u_0\leq b_0$ and $u_1\leq b_1$ in
$\Sigma$ such that $u\leq u_0\vee u_1$ is a \mjc. By
\HSj, there exists $i<2$ such that $u\leq p\vee u_i$ and either
$p\leq u_i\vee v,u\vee u_{1-i}$ or $p\leq u_0\vee v,u_1\vee v$. In
the first case,
   \[
   p\leq a\wedge(u_i\vee v)\wedge(u\vee u_{1-i})\wedge
   \Bigl(\bigl(u\wedge(p\vee u_i)\bigr)\vee v\Bigr)\leq e.
   \]
In the second case,
   \[
   p\leq a\wedge(u_0\vee v)\wedge(u_1\vee v)\wedge
   \Bigl(\bigl(u\wedge(p\vee u_i)\bigr)\vee v\Bigr)\leq e.
   \]
Since every element of $L$ is a join of elements of $\Sigma$, we
obtain that $d\leq e$. Since $e\leq d$ holds in any lattice, we
obtain that $d=e$.
\end{proof}

\begin{corollary}\label{C:HS2HSj}
The lattice $\Co(T)$ satisfies \HS, for every chain $(T,\utr)$.
\end{corollary}

\begin{proof}
We apply Lemma~\ref{L:HS2HSj} to $L=\Co(T)$ together with the
join-seed $\Sigma=\setm{\set{p}}{p\in T}$. Let $a$, $b$, $c$, $b_0$,
$b_1\in T$ such that $a\neq b$, the inequality
$\set{a}\leq\set{b}\vee\set{c}$ is minimal in $b$ (thus $a\neq c$),
and $\set{b}\leq\set{b_0}\vee\set{b_1}$. Since
$\Co(T,\utr)=\Co(T,\gtr)$, we may assume without loss of generality that
$c\tr a\tr b$.
Furthermore, there exists $i<2$ such that $b\utr b_i$,
whence $\set{b}\leq\set{a}\vee\set{b_i}$. Since $T$ is a chain,
either $b_{1-i}\utr a$ or $a\utr b_{1-i}$. In the first case,
$\set{a}\leq\set{b_i}\vee\set{c},\set{b}\vee\set{b_{1-i}}$. In the
second case, $\set{a}\leq\set{b_0}\vee\set{c},\set{b_1}\vee\set{c}$.

Hence $\Co(T)$ satisfies \HSj. By Lemma~\ref{L:HS2HSj}, $\Co(T)$
satisfies \HS.
\end{proof}

\section{The Transitivity Lemma}\label{S:Trans}

The main purpose of the present section is to prove the following
technical lemma, which provides a large supply of minimal coverings.

\begin{lemma}[The Transitivity Lemma]\label{L:Trans}
Let $L$ be a lattice satisfying the identities \HS, \Ud, \Bo, \El,
and \Pe, let $\Sigma$ be a join-seed of $L$, and let $a$, $b$, $c$,
$b_0$, $b_1\in\Sigma$ such that both $a\leq b\vee c$ and
$b\leq b_0\vee b_1$ are \mjc s. Then there exists $i<2$ such the
following statements hold:
\begin{enumerate}
\item the inequality $b\leq a\vee b_i$ holds, and both inequalities
$b\leq c\vee b_i$ and $a\leq c\vee b_i$ are \mjc s;

\item one of the following two statements holds:
\begin{itemize}
\item[(ii.1)] $a\leq b_i\vee c,b_{1-i}\vee b$ and, if $a\neq b_{1-i}$,
then the inequality $a\leq b_0\vee b_1$ is a \mjc;

\item[(ii.2)] $a\leq b_0\vee c,b_1\vee c$ and, if $a\neq b_{1-i}$,
then the inequality $a\leq b_{1-i}\vee c$ is a \mjc.
\end{itemize}
\end{enumerate}
\end{lemma}

The situation may be partly viewed on Figure~3.

\begin{proof}
It follows from Lemma~\ref{L:HS2HSj} that there exists $i<2$ such
that
   \begin{equation}\label{Eq:WhatMightHold}
   b\leq a\vee b_i\text{ and either }a\leq b_i\vee c,b_{1-i}\vee b
   \text{ or }a\leq b_0\vee c,b_1\vee c.
   \end{equation}
Since $b\leq b_i\vee c$ is a nontrivial join-cover and $\Sigma$ is a
join-seed of
$L$, there are $x\leq b_i$ and $c'\leq c$ in $\Sigma$ such that
$b\leq x\vee c'$
is a \mjc. By applying \Boj\ to the inequalities
$b\leq b_i\vee b_{1-i},x\vee c'$ and observing that
$b\nleq b_i=b_i\vee x$, we obtain that $b\leq b_{1-i}\vee x$,
whence, by the minimality assumption on $b_i$, $x=b_i$. By applying
Lemma~\ref{L:HS2HSj} to the \mjc s $a\leq b\vee c$ and
$b\leq b_i\vee c'$, we obtain that either $a\leq c'\vee c=c$,
\contr, or $a\leq b\vee c'$. By the minimality assumption on $c$,
the latter implies that $c=c'$. Hence we have proved the following:
   \begin{equation}\label{Eq:blcbimin}
   \text{the inequality }b\leq b_i\vee c\text{ is a \mjc}.
   \end{equation}
Now we shall proceed by proving the following statement:
   \begin{equation}\label{Eq:alcbimin}
   \text{the inequality }a\leq b_i\vee c\text{ is a \mjc}.
   \end{equation}
If $a\leq b_i$, then $b\leq a\vee b_i=b_i$, \contr; whence
$a\nleq b_i$. So $a\leq b_i\vee c$ is a nontrivial join-cover,
thus, since $\Sigma$ is a join-seed of $L$, there are $x\leq b_i$
and $c'\leq c$ in $\Sigma$ such that $a\leq x\vee c'$ is a \mjc. By
applying \Boj\ to the inequalities $a\leq b\vee c,x\vee c'$ and
observing that $a\nleq c=c\vee c'$, we obtain that $a\leq b\vee c'$,
whence, by the minimality assumption on $c$, we obtain that $c=c'$.

Now we apply Lemma~\ref{L:Pe2Pej} to the \mjc s
$a\leq b\vee c,x\vee c$ and the inequality $b\leq b_0\vee b_1$.
Thus either $b\leq a\vee x$ or there exists $j<2$ such that
$a\leq b_j\vee c$ and $b\leq a\vee b_j,b_j\vee x$.
Suppose that the second case holds. If $i\neq j$, then
$b\leq a\vee b_j\leq c\vee b_j$. But $b\leq c\vee b_i$ and
$b\leq b_i\vee b_j$, whence, by \Udj, either $b\leq b_0$ or
$b\leq b_1$ or $b\leq c$, \contr. Therefore, $i=j$ and
$b\leq x\vee b_i=b_i$, \contr.

Hence the first case holds, thus it follows from $a\leq x\vee c$ that
$b\leq x\vee c$ with $x\leq b_i$, thus, by~\eqref{Eq:blcbimin},
$x=b_i$. This completes the proof of~\eqref{Eq:alcbimin}, and thus
also the proof of (i).

Now let us establish the remaining \mjc s in (ii), under the additional
assumption that $a\neq b_{1-i}$. We have already seen that
$a\nleq b_i$. If $a\leq b_{1-i}$, then, since $b\leq a\vee b_i$ and by
the minimality assumption on $b_{1-i}$, we obtain that $a=b_{1-i}$,
\contr. Therefore, we have obtained the inequalities
   \begin{equation}\label{Eq:anlb0b1}
   a\nleq b_0\text{ and }a\nleq b_1.
   \end{equation}
Now we separate cases, according to \eqref{Eq:WhatMightHold}.
\smallskip

\noindent\textbf{Case 1.} $a\leq b_i\vee c,b\vee b_{1-i}$. {}From
the second inequality and $b\leq b_0\vee b_1$ it follows that
$a\leq b_0\vee b_1$. Thus, by \eqref{Eq:anlb0b1} and since $\Sigma$
is a join-seed of $L$, there are $x_0\leq b_0$ and $x_1\leq b_1$ in
$\Sigma$ such that $a\leq x_0\vee x_1$ is a \mjc. By applying \Boj\
to the inequalities $a\leq b_i\vee c$ (see \eqref{Eq:alcbimin}) and
$a\leq x_i\vee x_{1-i}$ and observing that $a\nleq b_i=x_i\vee b_i$,
we obtain the inequality $a\leq c\vee x_i$ with $x_i\leq b_i$,
thus, by \eqref{Eq:alcbimin}, $x_i=b_i$. On the other hand,
$b\leq a\vee b_i\leq b_i\vee x_{1-i}$ with $x_{1-i}\leq b_{1-i}$,
thus, by the minimality assumption on $b_{1-i}$, we obtain that
$x_{1-i}=b_{1-i}$. Therefore, we have proved the following
statement:
   \begin{equation}\label{Eq:alb0b1min}
   \text{the inequality }a\leq b_0\vee b_1\text{ is a \mjc}.
   \end{equation}

\smallskip
\noindent\textbf{Case 2.} $a\leq b_0\vee c,b_1\vee c$.
{}From \eqref{Eq:anlb0b1}, the inequalities $a\nleq c$ and
$a\leq b_{1-i}\vee c$, and the assumption that $\Sigma$ is a
join-seed of $L$, it follows that there are $x\leq b_{1-i}$ and
$c'\leq c$ in $\Sigma$ such that $a\leq x\vee c'$ is a \mjc. By
applying \Boj\ to the inequalities $a\leq b\vee c,x\vee c'$ and
observing that $a\nleq c=c\vee c'$, we obtain that $a\leq b\vee c'$,
whence, since $c'\leq c$ and by the minimality assumption on $c$,
we obtain that $c=c'$.

Suppose now that $x<b_{1-i}$.
Applying Lemma~\ref{L:Pe2Pej} to the join covers
$a\leq c\vee x,b\vee c$ and $b\leq b_0\vee b_1$, we obtain that
either $b\leq a\vee x$ or $b\leq a\vee b_j,x\vee b_j$, for some
$j<2$. In the first case, $b\leq a\vee x\leq c\vee x$. Since
$b\nleq b_{1-i}=b_{1-i}\vee x$, we obtain, by \Boj\ applied to the
inequalities $b\leq b_0\vee b_1,c\vee x$ that $b\leq b_i\vee x$, which
contradicts the assumption that the cover $b\leq b_0\vee b_1$ is a \mjc.
Hence the second case applies.  If $j\neq i$, then
$b\leq a\vee b_j\leq c\vee b_j$, while
$b\leq c\vee b_i$ and $b\leq b_i\vee b_j$, whence, by \Udj, either
$b\leq b_i$ or $b\leq b_j$ or $b\leq c$, \contr. Hence $j=i$ and
$b\leq x\vee b_i$ with $x<b_{1-i}$, which contradicts the minimality
assumption on $b_{1-i}$.
This completes the proof of the following statement:
   \begin{equation}\label{Eq:alb1-ic}
   \text{the inequality }a\leq b_{1-i}\vee c\text{ is a \mjc},
   \end{equation}
and thus the proof of (ii).
\end{proof}

In particular, in the context of Lemma~\ref{C:Trans}, it
follows from (i) that $a\DD b_i$ always holds. Moreover, if
$a\neq b_{1-i}$, then, by (ii), $a\DD b_{1-i}$ holds. Therefore, we
obtain the following remarkable corollary.

\begin{corollary}\label{C:Trans}
Let $L$ be a lattice satisfying the identities \HS, \Ud, \Bo, \El,
and \Pe, let $\Sigma$ be a join-seed of $L$. For any $a$, $b$,
$c\in\Sigma$, from $a\DD b\DD c$ and $a\neq c$ it follows that
$a\DD c$.
\end{corollary}

\section{The construction}\label{S:Constr}
In this section, we shall fix a complete, lower continuous, finitely
spatial lattice~$L$ satisfying \HS, \Ud, \Bo, \El, and \Pe. By
Lemma~\ref{L:HS2SDj}, $L$ is dually \IIdistr, thus, by
Lemma~\ref{L:JoinSeeds}, $\Sigma=\J(L)$ is a join-seed of $L$.

For every $a\in\J(L)$, we denote by $\set{A_a,B_a}$ the
Udav-Bond partition of $\rd{a}$ associated with $a$, as defined in
\cite[Section~5]{SeWe1}. We define a binary relation $\utra$ on
$\J_a(L)=\set{a}\cup\rd{a}$ by the following:
\begin{enumerate}
\item $x\utra a\utra y$ and $x\utra y$, for all
$(x,y)\in(A_a\cup\set{a})\times(B_a\cup\set{a})$;

\item $x\utra y$ if{f} $y\leq a\vee x$, for all $x$,
$y\in A_a$;

\item $x\utra y$ if{f} $x\leq a\vee y$, for all $x$,
$y\in B_a$.
\end{enumerate}

We also say that $x\tra y$ if{f} $x\utra y$ and $x\neq y$, for all
$x$, $y\in\J_a(L)$.

\begin{lemma}\label{L:utralo}
The relation $\utra$ is a total ordering of $\J_a(L)$, for any
$a\in\J(L)$.
\end{lemma}

\begin{proof}
It is trivial that $\utra$ is reflexive. Let $x$, $y$, $z\in\J_a(L)$
with $x\utra y$ and $y\utra z$, we prove that $x\utra z$. This is
obvious if either $a\in\set{x,y,z}$ or $(x,z)\in A_a\times B_a$, so
suppose otherwise. Then $x$ and $z$ belong to the same block of the
Udav-Bond partition associated with $a$, say,
$\set{x,z}\subseteq A_a$. Since $y\utra z$, $y$ belongs to $A_a$ as
well. Furthermore, $z\leq a\vee y\leq a\vee x$ and thus $x\utra z$. The
proof for $\set{x,z}\subseteq B_a$ is similar. This proves that
$\utra$ is transitive.

Let $x$, $y\in\J_a(L)$ such that $x\utra y\utra x$, we prove that
$x=y$. This is obvious if $a\in\set{x,y}$, so suppose that
$a\notin\set{x,y}$. Then $x$ and $y$ belong to the same block of the
Udav-Bond partition associated with $a$, say,
$\set{x,y}\subseteq A_a$. Pick $u\in B_a$. Then $a\vee x=a\vee y$,
but $a\leq u\vee x,u\vee y$, thus $u\vee x=u\vee y$, thus, by
the \jsdy\ of $L$ (see Lemma~\ref{L:HS2SDj}),
   \[
   a\leq u\vee x=u\vee y=u\vee(x\wedge y).
   \]
However, by Lemma~\ref{L:MinUB}, both inequalities
$a\leq u\vee x,u\vee y$ are \mjc s, thus $x=y$. Hence $\utra$ is
antisymmetric.

Now let $x$, $y\in\J_a(L)$, we prove that either $x\utra y$ or
$y\utra x$. This is obvious if either $a\in\set{x,y}$ or $x$ and $y$
belong to different blocks of the Udav-Bond partition associated with
$a$, so suppose otherwise, say, $\set{x,y}\subseteq B_a$. Pick
$u\in A_a$. By Lemma~\ref{L:MinUB}, both inequalities
$a\leq u\vee x,u\vee y$ are \mjc s, thus, by applying
Lemma~\ref{L:El2Elj} to the \mjc s
$a\leq u\vee x,u\vee y,u\vee y$, we obtain that either $x\leq a\vee y$
or $y\leq a\vee x$, thus either $x\utra y$ or $y\utra x$. The proof
for $\set{x,y}\subseteq A_a$ is similar. Hence $\utra$ is a total
ordering.
\end{proof}

For any $a\in\J(L)$, let $\gfa\colon L\to\PP(\J_a(L))$ be the
map defined by the rule
   \[
   \gfa(x)=\setm{b\in\J_a(L)}{b\leq x},\text{ for all }x\in L.
   \]

\begin{lemma}\label{L:gfaxConv}
The set $\gfa(x)$ is order-convex in $(\J_a(L),\utra)$, for any
$x\in L$.
\end{lemma}

\begin{proof}
Let $u$, $v$, $w\in\J_a(L)$ such that $u\tra w\tra v$ and $u$,
$v\leq x$, we prove that $w\leq x$. If $u\in\set{a}\cup A_a$ and
$v\in\set{a}\cup B_a$, then $a\leq u\vee v\leq x$, and then,
$w\in\set{a}\cup A_a$ implies that $w\leq a\vee u\leq x$, while
$w\in\set{a}\cup B_a$ implies that $w\leq a\vee v\leq x$.

Suppose now that $u$, $v\in A_a$. {}From $w\utra v$ it follows that
$w\in A_a$. Pick $t\in B_a$. By Lemma~\ref{L:MinUB}, all
inequalities $a\leq t\vee u,t\vee v,t\vee w$ are \mjc s; from
$a\vee v\leq a\vee w\leq a\vee u$ it follows that
$t\vee v\leq t\vee w\leq t\vee u$,
thus, by Lemma~\ref{L:IVP},
$w\leq u\vee v\leq x$. The argument is similar in case $u$,
$v\in B_a$.
\end{proof}

\begin{lemma}\label{L:LattHom}
The map $\gfa$ is a lattice homomorphism from $L$ to $\Co(\J_a(L))$,
and it preserves the existing bounds.
\end{lemma}

\begin{proof}
It is clear that $\gfa$ is a meet-homomorphism from $L$ to
$\Co(\J_a(L))$ and that it preserves the existing bounds. Let $x$,
$y\in L$, we prove that $\gfa(x\vee y)=\gfa(x)\vee\gfa(y)$. It
suffices to prove that $b\in\gfa(x\vee y)$ implies that
$b\in\gfa(x)\vee\gfa(y)$, for any $b\in\J_a(L)$. This is trivial if
$b\in\gfa(x)\cup\gfa(y)$, so suppose otherwise, that is,
$b\nleq x,y$. Since $b\leq x\vee y$ and $\J(L)$ is a join-seed of
$L$, there are $b_0\leq x$ and $b_1\leq y$ in $\J(L)$ such that the
inequality $b\leq b_0\vee b_1$ is a \mjc. {}From
Corollary~\ref{C:Trans} it follows that both $b_0$ and $b_1$ belong to
$\J_a(L)$. If $b=a$, then the pair $(b_0,b_1)$ belongs either to
$A_a\times B_a$ or $B_a\times A_a$. In the first case,
$b_0\utra a\utra b_1$, in the second case, $b_1\utra a\utra b_0$; in
both cases, $b=a\in\gfa(x)\vee\gfa(y)$.

Suppose from now on that $b\neq a$, say, $b\in B_a$. Pick $c\in\J(L)$
such that $a\leq b\vee c$ is a \mjc; observe that $c\in A_a$. So there
exists $i<2$ such that the statements (i), (ii) of Lemma~\ref{L:Trans}
hold.

{}From the fact that the inequality $a\leq b_i\vee c$ is a \mjc\ and
$c\in A_a$ it follows that $b_i\in B_a$. {}From the relations $b$,
$b_i\in B_a$, $b\neq b_i$, and $b\leq a\vee b_i$ it follows that
   \begin{equation}\label{Eq:btrabi}
   b\tra b_i.
   \end{equation}
If $a=b_{1-i}$, then, since $b\in B_a$, we obtain that
$b_{1-i}=a\tra b$, thus, by \eqref{Eq:btrabi},
$b\in\gfa(x)\vee\gfa(y)$. Suppose from now on that $a\neq b_{1-i}$.
If (ii.1) of Lemma~\ref{L:Trans} holds, then the inequality
$a\leq b_0\vee b_1$ is a \mjc\ with $b_i\in B_a$, thus
$b_{1-i}\in A_a$, thus $b_{1-i}\tra b$, which,
together with \eqref{Eq:btrabi}, implies that
$b\in\gfa(x)\vee\gfa(y)$. Suppose now that (ii.2) of
Lemma~\ref{L:Trans} holds. {}From the fact that $a\leq b_{1-i}\vee c$
is a \mjc\ and $c\in A_a$ it follows that $b_{1-i}\in B_a$. If
$b\utra b_{1-i}$, then, since $b$, $b_{1-i}\in B_a$, we obtain that
$b\leq a\vee b_{1-i}\leq c\vee b_{1-i}$, but $b\leq c\vee b_i$ and
$b\leq b_i\vee b_{1-i}$, whence, by \Udj, either $b\leq b_0$ or
$b\leq b_1$ or $b\leq c$, \contr. Hence $b\nutra b_{1-i}$, thus, by
Lemma~\ref{L:utralo}, $b_{1-i}\utra b$. Therefore, it follows again
from \eqref{Eq:btrabi} that $b\in\gfa(x)\vee\gfa(y)$.
\end{proof}

\section{The representation theorem}\label{S:RepThm}

\begin{notation}\label{No:ClassC}
Let $\CC$ denote the class of all lattices that can be embedded into
a direct product of the form $\prod_{i\in I}\Co(T_i)$,
where $\famm{T_i}{i\in I}$ is a family of chains.
\end{notation}

Our main theorem is the following.

\begin{theorem}\label{T:EmbThm}
For a lattice $L$, the following are equivalent:
\begin{enumerate}
\item $L$ belongs to $\CC$.

\item $L$ satisfies the identities \HS, \Ud, \Bo, \El, and \Pe.

\item There exists an embedding
$\gf\colon L\into\prod_{i\in I}\Co(T_i)$, for some family
$\famm{T_i}{i\in I}$ of chains, which preserves the existing bounds and
satisfies the following additional properties:
\begin{itemize}
\item[---] if $L$ is finite, then
$\sum_{i\in I}|T_i|\leq|\J(L)|^2$;

\item[---] if $L$ is subdirectly irreducible, then $I=\set{0}$,
$\gf$ is atom-preserving, and, if $L$ is finite, then
$|T_0|=|\J(L)|$.

\item[---] if $L$ is finite, atomistic, and subdirectly irreducible,
then $L\cong\Co(n)$, where $n=|\J(L)|$.
\end{itemize}
\end{enumerate}
\end{theorem}

\begin{proof}
(i)$\Rightarrow$(ii) We have seen in \cite{SeWe1} that $L$ satisfies
\Ud\ and \Bo. Moreover, it follows from Corollaries~\ref{C:El2Elj},
\ref{C:Pe2Pej}, and \ref{C:HS2HSj} that $L$ satisfies \El, \Pe, and
\HS.

(ii)$\Rightarrow$(iii) As in \cite{SeWe1,SeWe2}, we embed $L$ into the
filter lattice $\hL$ of $L$, partially ordered by reverse inclusion.
This embedding preserves the existing bounds and atoms.
We recall that $\hL$ is complete, lower continuous, and finitely
spatial. Let $\J_a(\hL)$ and $\gfa\colon\hL\to\Co(\J_a(\hL))$ be
defined as in Section~\ref{S:Constr}, and let
$\gya\colon L\to\Co(\J_a(\hL))$ be the restriction of $\gfa$ to $L$,
for any $a\in\J(\hL)$. Since every element of $\hL$ is a join of
elements of $\J(\hL)$, it follows from Lemma~\ref{L:LattHom} that
the map $\gy\colon L\to\prod_{a\in\J(\hL)}\Co(\J_a(\hL))$ that with
any $x\in L$ associates the family $\famm{\gya(x)}{a\in\J(\hL)}$ is
a lattice embedding; it obviously preserves the existing bounds.
In case $L$ is finite, we have $\hL=L$ and $|\J_a(L)|\leq|\J(L)|$, for
all $a\in\J(L)$; the cardinality bound follows immediately.

Suppose now that $L$ is subdirectly irreducible. Thus $\gya$ is an
embedding, for some $a\in\J(\hL)$; pick such an $a$. Every atom $x$ of
$L$ is also an atom of $\hL$, and $\gya(x)$ is nonempty, thus there
exists $b\in\J_a(\hL)$ below $x$, whence $x=b\in\J_a(\hL)$ and
$\gya(x)=\set{x}$, an atom of $\Co(\J_a(L))$.
Suppose now that $L$ is finite, thus $\hL=L$. For any $x\in\J(L)$,
if $x'$ denotes the join of all elements of $\J_a(L)$ below $x$,
then $\gya(x)=\gya(x')$, whence $x=x'$, thus, since $x$ is \jirr,
$x\in\J_a(L)$; therefore, $\J_a(L)=\J(L)$.

Now suppose, in addition, that $L$ is atomistic. Then
$\set{x}=\gya(x)$ belongs to the range of $\gya$, for any
$x\in\J(L)$, thus $\gya$ is surjective, hence it is an isomorphism
from~$L$ onto $\Co(\J(L),\utra)$.

(iii)$\Rightarrow$(i) is trivial.
\end{proof}

\begin{remark}\label{Re:EmbThm}
A finite, atomistic lattice $L$ in $\CC$ may not embed
atom-preservingly into any $\Co(P)$, thus \emph{a fortiori}
into any product of the form
$\prod_{i\in I}\Co(T_i)$ where the $T_i$-s are chains, as shows
\cite[Example~8.1]{SeWe1}. Also, a finite, atomistic, subdirectly
irreducible lattice in $\SUB$ may not be of the form $\Co(P)$, see
\cite[Example~8.2]{SeWe1}.
\end{remark}

\begin{corollary}\label{C:EmbThm}
The class $\CC$ is a finitely based variety of lattices. In
particular, $\CC$ is closed under homomorphic images.
\end{corollary}

This result solves positively Problem~3 in \cite{SeWe1}.

\section{The class $\SUB(n)$, for $n\geq 0$}\label{S:SUB(n)}

We start with the following lemma.

\begin{lemma}\label{L:Min2St}
Let $L$ be a complete, lower continuous, finitely spatial lattice in
$\CC$, let $a\in\J(L)$. Let $x$, $y$, $u\in\rd{a}$ such that
$a\leq u\vee x,u\vee y$. If $x\leq a\vee y$, then the inequality
$x\leq u\vee y$ is a \mjc.
\end{lemma}

\begin{proof}
{}From Lemma~\ref{L:MinUB} it follows that both inequalities
$a\leq u\vee x,u\vee y$ are \mjc s.
Since $x\leq a\vee y$ and $a\leq u\vee y$, we obtain that
$x\leq u\vee y$. {}From Lemma~\ref{L:rdaAnt} it follows that
$x\nleq u,y$. Since $x\leq u\vee y$ and $\J(L)$ is a join-seed of
$L$, there are $u'\leq u$ and $y'\leq y$ in $\J(L)$ such that the
inequality $x\leq u'\vee y'$ is a \mjc. So
$a\leq x\vee u\leq y'\vee u$ with $y'\leq y$, thus, by the minimality
of $y$ in $a\leq y\vee u$, we obtain that $y'=y$. If $u'=a$, then
$a\leq u$, \contr; whence $u'\neq a$; but $a\DD x\DD u'$, whence, by
Corollary~\ref{C:Trans}, $a\DD u'$. But $u'\leq u$ and $a\DD u$,
whence, by Lemma~\ref{L:rdaAnt}, $u'=u$.
\end{proof}

Now we are able to relate chains in the $\J_a(L)$-s and Stirlitz
tracks.

\begin{corollary}\label{C:Lh2St}
Let $L$ be a complete, lower continuous, finitely spatial lattice in
$\CC$, let $a\in\J(L)$, let $n$ be a natural number, let $u$, $x_0$,
\dots, $x_n\in\J_a(L)$ with $x_0\tra x_1\tra\cdots\tra x_n$. Denote by
$\set{A_a,B_a}$ the Udav-Bond partition of $\rd{a}$ associated with
$a$. Then the following statements hold:
\begin{enumerate}
\item If $u\in A_a$ and $x_0$, \dots, $x_n\in\set{a}\cup B_a$, then
$(\famm{x_i}{0\leq i\leq n},\famm{u}{1\leq i\leq n})$ is a Stirlitz
track.

\item If $u\in B_a$ and $x_0$, \dots, $x_n\in\set{a}\cup A_a$, then
$(\famm{x_{n-i}}{0\leq i\leq n},\famm{u}{1\leq i\leq n})$ is a
Stirlitz track.
\end{enumerate}
\end{corollary}

\begin{proof}
(i) It follows from Lemma~\ref{L:Min2St} that the inequality
$x_i\leq u\vee x_{i+1}$ is a \mjc, for any $i\in\fso{n-1}$; the
conclusion follows. The proof for (ii) is similar.
\end{proof}

We recall, see \cite{SeWe2}, that for any positive integer $n$, the
class $\SUB_n$ of all lattices that can be embedded into some $\Co(P)$
where $P$ is a poset of length at most $n$ is a finitely based
variety, defined by the identities \St, \Ud, \Bo, together with new
identities \Ht{n} and \Ht{k,n+1-k} for $1\leq k\leq n$.

\begin{notation}
For a natural number $n$, let $\SUB(n)$ denote the class of all
lattices that can be embedded into a power of $\Co(n)$.
\end{notation}

Of course, $\SUB(0)$ is the trivial variety while $\SUB(1)=\SUB(2)$ is
the class of all distributive lattices. Now we obtain the main result
of this section.

\begin{theorem}\label{T:SUB(n)}
Let $n$ be a positive integer.
The class $\SUB(n+1)$ is a finitely generated variety, defined by the
identities \HS, \Ud, \Bo, \El, \Pe, and \Ht{k,n+1-k} for
$1\leq k\leq n$. Furthermore, $\SUB(n+1)=\CC\cap\SUB_n$.
\end{theorem}

\begin{proof}
Since the $(n+1)$-element chain belongs to $\SUB_n$,
the containment\linebreak
$\SUB(n+1)\subseteq\CC\cap\SUB_n$ is
obvious. Furthermore, by the results of~\cite{SeWe2} and
Theorem~\ref{T:EmbThm}, every lattice in
$\CC\cap\SUB_n$ satisfies the identities
\HS, \Ud, \Bo, \El, \Pe, and \Ht{k,n+1-k} for $1\leq k\leq n$.

Now let $L$ be a lattice satisfying the identities
\HS, \Ud, \Bo, \El, \Pe, and \Ht{k,n+1-k} for $1\leq k\leq n$, we
prove that $L$ belongs to $\SUB(n+1)$. By embedding~$L$ into its
filter lattice, we see that it suffices to consider the case where
$L$ is complete, lower continuous, and finitely spatial. By
Theorem~\ref{T:EmbThm}, $L$ belongs to $\CC$. In order to conclude
the proof, it suffices to establish that $\J_a(L)$ has at most $n+1$
elements, for any $a\in\J(L)$. If this is not the case, then
both blocks $A_a$ and $B_a$ of the Udav-Bond partition of $\rd{a}$
associated with $a$ are nonempty, and $\J_a(L)$ has a chain of the
form
   \[
   x_k\tra\cdots\tra x_1\tra x_0=a=y_0\tra y_1\tra\cdots\tra y_l,
   \]
where $k$ and $l$ are positive integers with $k+l=n+1$. Define pairs
$\sigma$ and $\tau$ by
   \begin{align*}
   \sigma&=(\famm{x_i}{0\leq i\leq k},\famm{y_1}{1\leq i\leq k}),\\
   \tau&=(\famm{y_j}{0\leq j\leq l},\famm{x_1}{1\leq j\leq l}).
   \end{align*}
It follows from Corollary~\ref{C:Lh2St} that both $\sigma$ and $\tau$
are Stirlitz tracks, but $a\leq x_1\vee y_1$, thus the pair
$(\sigma,\tau)$ is a bi-Stirlitz track (see \cite{SeWe2}) of index
$(k,l)$ with $k+l=n+1$, which contradicts the fact that $L$ satisfies
the identity \Ht{k,l}, see \cite[Proposition~6.2]{SeWe2}.

In particular, we have proved that $\SUB(n+1)$ is a variety.
Of course, it is generated by the single finite lattice $\Co(n+1)$.
\end{proof}

Since the construction underlying Theorem~\ref{T:SUB(n)} is the same
as the one underlying Theorem~\ref{T:EmbThm}, the corresponding
additional information is preserved. For example, any member $L$ of
$\SUB(n+1)$ has an embedding into a power of $\Co(n+1)$ which
preserves the zero if it exists; furthermore, if $L$ is subdirectly
irreducible, then this embedding preserves atoms.

\begin{theorem}\label{T:CLocFin}
The variety $\CC$ is locally finite.
\end{theorem}

\begin{proof}
For a lattice $L$, let $\Csub(L)$ denote the lattice of all convex
sublattices of~$L$, ordered by inclusion. For a variety $\VV$ of
lattices, let $\Csub(\VV)$ denote the variety generated by all
lattices of the form $\Csub(L)$, for $L\in\VV$. For a chain $T$, the
equality $\Co(T)=\Csub(T)$ obviously holds, whence $\CC$ is a
subvariety of $\Csub(\bD)$, where $\bD$ denotes the variety of all
distributive lattices. It is proved in V.~Slav{\'\i}k~\cite{Slav} that
$\Csub(\bD)$ is locally finite, therefore, the smaller variety $\CC$
is also locally finite.
\end{proof}

\begin{corollary}\label{C:CLocFin}
The variety $\CC$ is generated by $\Co(\omega)$, where $\omega$ denotes
the chain of natural numbers.
\end{corollary}

If, for a poset $P$, we denote by $\SUB(P)$ the variety generated by
$\Co(P)$, we obtain the `equation' $\CC=\SUB(\omega)$.

\begin{corollary}\label{C:QVjoin}
The variety $\CC$ is the quasivariety join of all varieties
$\SUB(n)$, where $1\leq n<\omega$.
\end{corollary}

\begin{proof}
Let $\QQ$ be any quasivariety containing $\SUB(n)$, for every positive
integer $n$. Every finite lattice $L$ in $\CC$ embeds into a finite
power of some $\Co(n)$, thus it belongs to $\QQ$. By
Theorem~\ref{T:CLocFin}, it follows that $\QQ$ contains $\CC$.
\end{proof}

\section{Weak Stirlitz tracks in lattices of convex
subsets of chains}\label{S:StSubC}

\begin{definition}\label{D:wSt}
Let $L$ be a lattice, let $m$, $n$ be positive integers.
\begin{enumerate}
\item A \emph{weak Stirlitz track} of length $n$ of $L$ is a pair
$\sigma=(\famm{x_i}{0\leq i\leq n},x)$, where $x$, $x_i$ (for
$0\leq i\leq n$) are elements of $L$, and the following relations
hold:
   \begin{itemize}
\item[(1)] $x_0\neq(x_0\wedge x_1)\vee(x_0\wedge x)$;

\item[(2)] $x_k\leq x_{k+1}\vee x$, for all $k\in\fso{n-1}$;

\item[(3)] $x_{k-1}\nleq(x_k\wedge x_{k+1})\vee x$, for all
$k\in\fsi{n-1}$.
   \end{itemize}
\item A \emph{weak bi-Stirlitz track} of index $(m,n)$ of $L$ is a
pair $(\sigma,\tau)$, where $\sigma=(\famm{x_i}{0\leq i\leq m},x)$
and $\tau=(\famm{y_j}{0\leq j\leq n},y)$ are both weak Stirlitz
tracks such that $x_0=y_0\leq x_1\vee y_1$ while
$x_0\neq(x_0\wedge x_1)\vee(x_0\wedge y_1)$.
\end{enumerate}
\end{definition}

For a Stirlitz track
$\sigma=(\famm{x_i}{0\leq i\leq n},\famm{x'_i}{1\leq i\leq n})$, we
put \linebreak $\ol{\sigma}=(\famm{x_i}{0\leq i\leq n},x'_1)$, and
$\tilde{\sigma}=\famm{x_i}{0\leq i\leq n}$, the \emph{trace}
of $\sigma$ and of $\ol{\sigma}$. The trace of a (weak) bi-Stirlitz track
$(\sigma,\tau)$ is the pair $(\tilde{\sigma},\tilde{\tau})$.

\begin{lemma}\label{L:wSt}
Let $L$ be a lattice in $\SUB$. Then the following statements hold:
\begin{enumerate}
\item The pair $\ol{\sigma}$ is a weak Stirlitz track of
$L$, for every Stirlitz track $\sigma$ of $L$.

\item The pair $(\ol{\sigma},\ol{\tau})$ is a weak
bi-Stirlitz track of $L$, for every bi-Stirlitz track
$(\sigma,\tau)$ of $L$.
\end{enumerate}
\end{lemma}

\begin{proof}
(i) Let
$\sigma=(\famm{x_i}{0\leq i\leq n},\famm{x'_i}{1\leq i\leq n})$ be
a Stirlitz track of $L$. We put $x=x'_1$ and we verify (1)--(3) of
Definition~\ref{D:wSt}(i).

The inequality (1) is trivial, while the inequality (2) follows
from \cite[Lemma~5.6]{SeWe1}. Suppose that
$x_{k-1}\leq(x_k\wedge x_{k+1})\vee x'_1$. If $x_{k-1}\leq x'_1$,
then, again by \cite[Lemma~5.6]{SeWe1},
$x_0\leq x_{k-1}\vee x'_1=x'_1$, \contr; whence
$x_{k-1}\nleq x'_1$. Furthermore, $x_{k-1}\nleq x_k$, thus
$x_{k-1}\nleq x_k\wedge x_{k+1}$. By applying \Bo\ to the
inequalities
$x_{k-1}\leq x_k\vee x'_k,(x_k\wedge x_{k+1})\vee x'_1$ and
observing that $x_{k-1}\nleq x_k=x_k\vee(x_k\wedge x_{k+1})$, we
obtain that $x_{k-1}\leq(x_k\wedge x_{k+1})\vee x'_k$ with
$x_k\wedge x_{k+1}<x_k$, which contradicts the minimality assumption on
$x_k$.

(ii) follows immediately from (i) and the \jirry\ of $x_0$.
\end{proof}

For subsets $X$ and $Y$ of a chain $(T,\utr)$, let
$X\tr Y$ and $X\bel Y$ be the following statements:
   \begin{align*}
   X\tr Y&\rightleftharpoons x\tr y,
   \text{ for all }(x,y)\in X\times Y,\\
   X\bel Y&\rightleftharpoons
   \forall x\in X,\ \exists y\in Y\text{ such that }x\utr y.
   \end{align*}
Of course, the equivalence
   \[
   X\tr Y\Longleftrightarrow(X\bel Y\text{ and }X\cap Y=\es)
   \]
holds, for all nonempty $X$, $Y\in\Co(T)$.

\begin{lemma}\label{L:DisjInt}
Let $(T,\utr)$ be a chain, let $L$ be a sublattice of $\Co(T)$, let
$n$ be a positive integer. For any weak Stirlitz track
$(\famm{X_i}{0\leq i\leq n},X)$ of~$L$, either the following statement
or its dual holds:
   \[
   X\bel X_0\bel X_1\text{ and }X\tr X_1\tr X_2\tr\cdots\tr X_n.
   \]
\end{lemma}

\begin{proof}
If $X_1\cap X\neq\es$, then
$Y\cap(X_1\vee X)=(Y\cap X_1)\vee(Y\cap X)$ for any
$Y\in\Co(T)$, \contr\ for $Y=X_0$. Hence $X_1\cap X=\es$. It follows
that either $X\tr X_1$ or $X_1\tr X$, say, $X\tr X_1$. Since
$X_0\leq X_1\vee X$ is a nontrivial join-cover, we obtain that
$X\bel X_0\bel X_1$.

Now we prove, by induction on $k$, that the statement
   \begin{equation}\label{Eq:belIncr}
   X\tr X_1\tr X_2\tr\cdots\tr X_k
   \end{equation}
holds, for any $k\in\fsi{n}$. For $k=1$ this is already verified.
Suppose having established \eqref{Eq:belIncr} at step $k$, with
$1\leq k<n$.

Suppose that $X_k\not\bel X_{k+1}$, that is, there exists $x\in X_k$
such that $y\tr x$ holds for any $y\in X_{k+1}$. By the induction
hypothesis, this also holds for any $y\in X$, thus
$x\notin X\vee X_{k+1}$, which contradicts the assumption that
$X_k\subseteq X\vee X_{k+1}$. Hence we have proved the relation
   \begin{equation}\label{Eq:Xkk+1bel}
   X_k\bel X_{k+1}.
   \end{equation}
Suppose that $X_k\cap X_{k+1}\neq\es$. Since
$X\bel X_{k-1}\bel X_k\bel X_{k+1}$, we obtain that
$X_{k-1}\subseteq(X_k\cap X_{k+1})\vee X$, \contr. Hence we have
established the relation
   \begin{equation}\label{Eq:Xkk+1es}
   X_k\cap X_{k+1}=\es.
   \end{equation}
{}From \eqref{Eq:Xkk+1bel} and \eqref{Eq:Xkk+1es} it follows that
$X_k\tr X_{k+1}$, which completes the induction step for
\eqref{Eq:belIncr}. For $k=n$, we obtain the conclusion of the lemma.
\end{proof}

\begin{lemma}\label{L:DisjInt2}
Let $(T,\utr)$ be a chain, let $L$ be a sublattice of $\Co(T)$, let
$m$ and $n$ be positive integers, let $(\sigma,\tau)$ be a weak
bi-Stirlitz track of $\Co(T)$ of index $(m,n)$, with
   \begin{align*}
   \sigma&=(\famm{X_i}{0\leq i\leq m},X),\\
   \tau&=(\famm{Y_j}{0\leq j\leq n},Y).
   \end{align*}
Then $X_1\vee Y_1\neq X_1\cup Y_1$, and, putting $Z=X_0=Y_0$,
either the following statement or its dual holds:
   \begin{align*}
   X_m\tr\cdots\tr X_1\tr Y_1\tr\cdots\tr Y_n\text{ and }
   X_1\bel Z\bel Y_1.
   \end{align*}
Furthermore, $Z$ does not meet simultaneously $X_1$ and $Y_1$,
and $\Co(m+n)$ embeds into~$L$.
\end{lemma}

\begin{proof}
It follows from
Lemma~\ref{L:DisjInt} that we may assume, without loss of
generality, that the following statement holds:
   \begin{equation}\label{Eq:belDecrX}
   X_m\tr\cdots\tr X_1\tr X\text{ and }X_1\bel X_0\bel X.
   \end{equation}
Suppose that $X_0\bel X_1$. Since $X_1\bel X_0$ and
$X_0\not\subseteq X_1$, $X_1$ is a proper final segment of $X_0$.
Thus, from $X_1\tr X$ it follows that $X_0\tr X$, but
$X_0\subseteq X_1\vee X$, whence $X_0\subseteq X_1$, \contr. Hence
we have established the relation
   \begin{equation}\label{Eq:X0notbelX1}
   X_0\not\bel X_1.
   \end{equation}
Now suppose that $Y_1\bel Y_0$ and $Y_1\tr Y$. As in the paragraph
above, we obtain that $Y_0\not\bel Y_1$. By \eqref{Eq:X0notbelX1} and
since $X_0=Y_0=Z$, there exists $z\in Z$ such that $y\tr z$ for any
$y\in X_1\cup Y_1$, which contradicts the fact that
$Z\subseteq X_1\vee Y_1$. Therefore, by Lemma~\ref{L:DisjInt}, the
following statement holds:
   \begin{equation}\label{Eq:belIncrY}
   Y\tr Y_1\tr\cdots\tr Y_n\text{ and }Y\bel Y_0\bel Y_1.
   \end{equation}
If $X_1\vee Y_1=X_1\cup Y_1$, then $Z=(Z\cap X_1)\vee(Z\cap Y_1)$,
\contr. Hence $X_1\vee Y_1\neq X_1\cup Y_1$, in particular,
$X_1\cap Y_1=\es$. This, together with
\eqref{Eq:belDecrX} and \eqref{Eq:belIncrY}, establishes the statement
   \[
   X_m\tr\cdots\tr X_1\tr Y_1\tr\cdots\tr Y_n.
   \]
Furthermore, if $Z$ meets both
$X_1$ and $Y_1$, then $Z=(Z\cap X_1)\vee(Z\cap Y_1)$, \contr.

In particular, sending $\set{i}$ to $X_{m-i}$ for $0\leq i<m$ and to
$Y_{i-m+1}$ for $m\leq i<m+n$ defines a lattice embedding from
$\Co(m+n)$ into $L$.
\end{proof}

\section{Subvarieties of $\CC$}\label{S:SubvC}

\begin{notation}\label{No:Lmn}
For positive integers $m$ and $n$, we set
   \[
   L_{m,n}=\setm{X\in\Co(m+n+1)}{m\in X\Rightarrow m-1\in X},
   \]
and we put $c_m=\set{m-1,m}$. Observe that $c_m\in\J(L_{m,n})$.
\end{notation}

The lattices $L_{m,n}$, for $m+n\leq 4$, are diagrammed on
Figure~4, together with $\Co(3)$ and $\Co(4)$.

\begin{definition}\label{D:CanWStrack}
For positive integers $m$ and $n$, the \emph{canonical bi-Stirlitz
track} of~$L_{m,n}$ is defined as $(\sigma_0,\tau_0)$, where we put
   \begin{align*}
   \sigma_0&=(\seq{c_m,\set{m-1},\dots,\set{0}},
   \famm{\set{m+n}}{1\leq i\leq m}),\\
   \tau_0&=\seq{c_m,\set{m+1},\dots,\set{m+n}},
   \famm{\set{0}}{1\leq j\leq n}).
   \end{align*}
\end{definition}
We observe that the relation $\set{m-1}<c_m$ (between entries of
$(\sigma_0,\tau_0)$) holds.

\begin{lemma}\label{L:BasicLmn}
Let $m$ and $n$ be positive integers. Then the following statements
hold:
\begin{enumerate}
\item $L_{m,n}$ is a $\seq{0,1}$-sublattice of $\Co(m+n+1)$.

\item The \jirr\ elements of $L_{m,n}$
are the singletons $\set{i}$, where $0\leq i\leq m+n$ and
$i\neq m$, together with $c_m$.

\item $L_{m,n}$ is subdirectly irreducible, with monolith
\pup{smallest nonzero congruence} $\Theta(\set{m-1},c_m)$.

\item All weak bi-Stirlitz tracks $\gamma$ of $L_{m,n}$ with
index $(m',n')$ such that $m'+n'=m+n$ have trace either
$(\tilde{\sigma}_0,\tilde{\tau}_0)$ or
$(\tilde{\tau}_0,\tilde{\sigma}_0)$ \pup{see
Definition~\textup{\ref{D:CanWStrack}}}.
\end{enumerate}
\end{lemma}

\begin{proof}
(i)--(iii) are straightforward. The result of (iv) follows easily
from Lemma~\ref{L:DisjInt2}.
\end{proof}

The proof of the following lemma is straightforward.

\begin{lemma}\label{L:betaSt}
Let $K$ and $L$ be lattices, let $f\colon K\onto L$ be a
lower bounded, surjective lattice homomorphism, let
$\beta\colon L\into K$ be the \jh\ defined by
$\beta(x)=\min f^{-1}\set{x}$, for all $x\in L$. Then the following
statements hold:
\begin{enumerate}
\item The image under $\beta$ of $\J(L)$ is contained in $\J(K)$.

\item The image under $\beta$ of any \mjc\ of $L$ is a \mjc\ of $K$.

\item The image under $\beta$ of any Stirlitz track \pup{resp.,
bi-Stirlitz track} of $L$ is a Stirlitz track \pup{resp.,
bi-Stirlitz track} of $K$.
\end{enumerate}
\end{lemma}

Now we can classify all finite subdirectly irreducible members of
$\CC$.

\begin{theorem}\label{T:ClassSI}\hfill
\begin{enumerate}
\item Let $L$ be a finite subdirectly irreducible lattice in $\CC$, let
$n$ be a positive integer. Then either $\Co(n)$ embeds into $L$ or $L$
embeds into $\Co(n)$.

\item Let $\VV$ be a subvariety of $\CC$. Then either
$\SUB(n)\subseteq\VV$ or $\VV\subseteq\SUB(n)$, for every positive
integer $n$.

\item The only finite subdirectly irreducible members of $\CC$ are
the $\Co(n)$, for $n>0$, and the $L_{m,n}$, for $m$, $n>0$.
\end{enumerate}
\end{theorem}

\begin{proof}
(i) Suppose that $L$ does not embed into $\Co(n)$. Since $L$ is
subdirectly irreducible, it does not belong to $\SUB(n)$, thus, by
Theorem~\ref{T:SUB(n)}, it does not belong to $\SUB_{n-1}$. Hence, there
exists $k\in\fsi{n-1}$ such that $L$ does not satisfy the identity
\Ht{k,n-k}, see~\cite{SeWe2}. Since $L$ is finite, it follows from
\cite[Proposition~6.2]{SeWe2} that $L$ has a bi-Stirlitz track of index
$(k,n-k)$, thus, by Lemma~\ref{L:DisjInt2}, $\Co(n)$ embeds into $L$.

(ii) Suppose that $\SUB(n)$ is not contained in $\VV$, that is,
$\Co(n)\notin\VV$. We prove that any lattice $L\in\VV$ belongs to
$\SUB(n)$. Since $\CC$ is locally finite (Theorem~\ref{T:CLocFin}),
it suffices to consider the case where $L$ is finite, hence it
suffices to consider the case where $L$ is finite and subdirectly
irreducible. {}From $\Co(n)\notin\VV$ it follows that $\Co(n)$ does not
embed into $L$, thus, by (i), $L$ embeds into $\Co(n)$, thus it belongs
to $\SUB(n)$.

(iii) Let $L$ be a finite subdirectly irreducible member of $\CC$.
Suppose that $L$ is nondistributive. There exists a largest
integer $n\geq 2$ such that $\Co(n)$ embeds into $L$. By (i), $L$
embeds into $\Co(n+1)$. Suppose that $L$ is not isomorphic to
$\Co(n)$. Since $L$ is subdirectly irreducible,
$L\notin\SUB(n)$, thus, as in the proof of (i), there
are $k$, $l>0$ such that $k+l=n$ and $L$ does not satisfy
\Ht{k,l}, hence $L$ has a bi-Stirlitz track $(\sigma,\tau)$ of
index $(k,l)$, with, say,
   \begin{align*}
   \sigma&=(\famm{X_i}{0\leq i\leq k},\famm{X'_i}{1\leq i\leq k}),\\
   \tau&=(\famm{Y_j}{0\leq j\leq l},\famm{Y'_j}{1\leq j\leq l}).
   \end{align*}
Put $Z=X_0=Y_0$. It follows from Lemma~\ref{L:DisjInt2} that, up to
possibly reversing the ordering of $n+1$ or exchanging $\sigma$
and $\tau$,
   \begin{gather}
   X_k\tr\cdots\tr X_1\tr Y_1\tr\cdots\tr Y_l\text{ and }
   X_1\bel Z\bel Y_1,
   \label{Eq:StTrNor}\\
   X_1\vee Y_1\neq X_1\cup Y_1\text{ and }Z\cap Y_1=\es.
   \label{Eq:ZX1Y1disj}
   \end{gather}
Since~$L$ has at most $n+1$ \jirr\ elements,
these elements are exactly the $X_i$-s, for $1\leq i\leq k$, the
$Y_j$-s, for $1\leq j\leq l$, and $Z$. Furthermore, it follows
from \eqref{Eq:StTrNor} and \eqref{Eq:ZX1Y1disj} that
$X_i=\set{k-i}$ for $1\leq i\leq k$,
$Y_j=\set{k+j}$ for $1\leq j\leq l$, and $Z$ is either equal to
$\set{k}$ or to $\set{k-1,k}$. In the first case, $L\cong\Co(n+1)$,
\contr, thus the second case applies. But then, $L\cong L_{k,l}$.
\end{proof}

\begin{remark}\label{Rk:InfSI}
There exists a proper class of infinite subdirectly irreducible
lattices in $\CC$, for example, all lattices of the form $\Co(T)$ where
$T$ is an infinite chain. However, each of those lattices generates the
variety $\CC$.
\end{remark}

As the union of the $\SUB(n)$, for $1\leq n<\omega$, generates
$\CC$, we obtain the following corollary.

\begin{corollary}\label{C:ClassSI}
Every proper subvariety of $\CC$ is finitely generated.
\end{corollary}

For a lattice $L$, let $\Var(L)$ denote the lattice variety
generated by $L$.

\begin{proposition}\label{P:MutInc}
Let $(m,n)$ and $(m',n')$ be pairs of positive integers such that
$m+n=m'+n'$. If $L_{m,n}$ belongs to $\Var(L_{m',n'})$, then
$(m,n)=(m',n')$.
\end{proposition}

\begin{proof}
By J\'onsson's Lemma (see B. J\'onsson \cite{Jons67} or P. Jipsen and
H. Rose \cite{JiRo}), there are a sublattice $L$ of
$L_{m',n'}$ and a congruence $\theta$ of $L$ such that
$L_{m,n}\cong L/{\theta}$.
The canonical bi-Stirlitz track $(\sigma_0,\tau_0)$ of index $(m,n)$
of $L_{m,n}\cong L/{\theta}$ can be, by Lemma~\ref{L:betaSt}, lifted
to a bi-Stirlitz track $(\sigma,\tau)$ of index $(m,n)$ of $L$, say,
   \begin{align*}
   \sigma&=(\famm{x_i}{0\leq i\leq m},
   \famm{x'}{1\leq i\leq m}),&&\text{for some }x'\in\J(L),\\
   \tau&=(\famm{y_j}{0\leq j\leq n},\famm{y'}{1\leq j\leq n}),
   &&\text{for some }y'\in\J(L),
   \end{align*}
with the additional property
   \begin{equation}\label{Eq:x1<x0}
   x_1<x_0
   \end{equation}
(because $(\sigma_0,\tau_0)$ has this property and
the map $\beta$ of Lemma~\ref{L:betaSt} is an order-embedding).
By Lemma~\ref{L:wSt}, $(\ol{\sigma},\ol{\tau})$ is a weak
bi-Stirlitz track of $L$, thus of $L_{m',n'}$, of index $(m,n)$, thus,
by Lemma~\ref{L:BasicLmn}(iv), its trace is either
$(\tilde{\sigma}_0,\tilde{\tau}_0)$ or
$(\tilde{\tau}_0,\tilde{\sigma}_0)$. But by \eqref{Eq:x1<x0}, only the
first case is possible, whence $(m,n)=(m',n')$.
\end{proof}

\begin{corollary}\label{C:subSUB(n)}
For any integer $n\geq 2$, the lattice $\boldsymbol{B}_n$ of
all lattice varieties $\VV$ such that
$\SUB(n)\subseteq\VV\subset\SUB(n+1)$ is isomorphic to $\two^{n-1}$.
\end{corollary}

\begin{proof}
It follows from Theorem~\ref{T:ClassSI} that the \jirr\ elements of
$\boldsymbol{B}_n$ are exactly the varieties $\Var(L_{k,l})$, where
$k$, $l>0$ and $k+l=n$. Furthermore, by Proposition~\ref{P:MutInc},
these varieties are mutually incomparable, hence they are atoms of
$\boldsymbol{B}_n$. Since $\boldsymbol{B}_n$ is finite
distributive, it is Boolean with $n-1$ atoms.
\end{proof}

The results of this section describe completely the lattice of
all subvarieties of $\CC$. This lattice is countable. Its bottom
is diagrammed on the left half of Figure~4. We use standard
notation, for example, $\mathbf{N}_5$ denotes the variety generated
by the pentagon, $\mathbf{L}_{1,2}$ denotes the variety generated by
$L_{1,2}$, and so on. The right half of Figure~4 represents small
subdirectly irreducible members of $\CC$.

\begin{figure}[htb]
\includegraphics{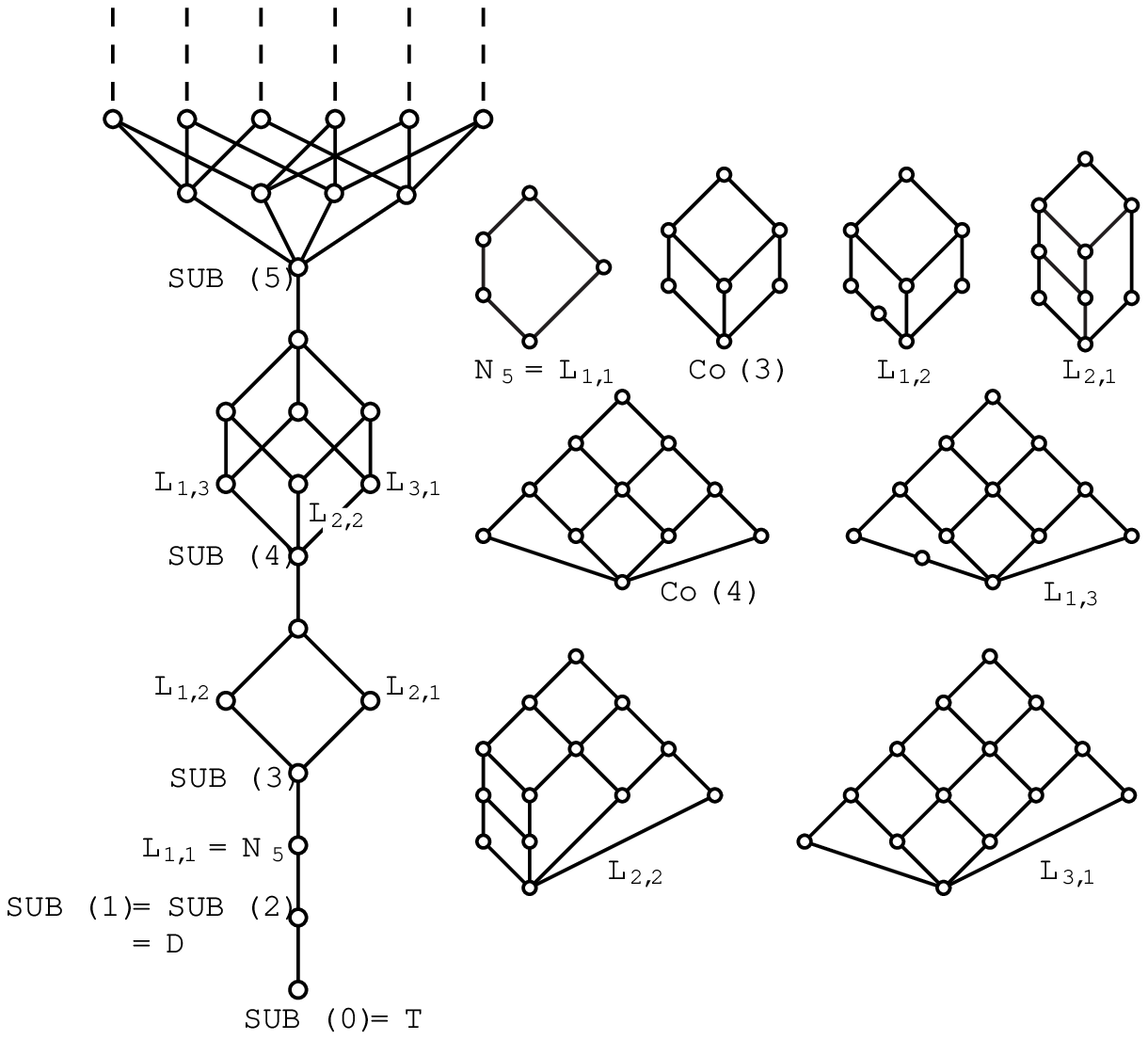}
\caption{Small subvarieties of $\CC$}
\end{figure}

\section{Projective members of $\CC$}\label{S:Proj}

\begin{notation}\label{No:Lambdamn}
Let $m$, $n>0$. We define lattice-theoretical statements
$\Lambda_n(x_0,\dots,x_{n-1})$ and $\Lambda_{m,n}(x_0,x_1,\dots,x_{m+n})$
as follows:
   \begin{align*}
   \Lambda_n(x_0,\dots,x_{n-1})\rightleftharpoons&
   \ x_k\leq x_i\vee x_j\text{ if }0\leq i<k<j<n\\
   &\text{ and } x_i\wedge x_j=x_0\wedge x_1\text{ for }i\neq j;\\
   \Lambda_{m,n}(x_0,x_1,\dots,x_{m+n})\rightleftharpoons&
   \ x_k\leq x_i\vee x_j\text{ if }0\leq i<k<j\leq m+n,\\
   &\ x_{m-1}\leq x_m,\\
   &\text{ and } x_i\wedge x_j=x_0\wedge x_2\text{ for }i\neq j\text{ and }
   \set{i,j}\neq\set{m-1,m}.
   \end{align*}
\end{notation}

We leave to the reader the easy proof of the following lemma.

\begin{lemma}\label{L:ChInt}
Let $(T,\utr)$ be a finite chain, let $n>0$, let $A_0$, \dots,
$A_{n-1}$ be pairwise disjoint elements of $\Co(T)$ such that
$A_k\subseteq A_i\vee A_j$, for $0\leq i<k<j<n$. Then either the
following statement or its dual holds:
\begin{quote}
There are elements $x_i$, $y_i$ $(i<n)$ of $T$ such that $A_i=[x_i,y_i)$,
for all $i<n$, and
   \[
   x_0\utr y_0\utr x_1\utr\cdots\utr x_{n-1}\utr y_{n-1}.
   \]
\end{quote}
\end{lemma}

The following lemma is the key to all projectivity results of the
present section.

\begin{lemma}\label{L:Lambdamn}
Let $L\in\CC$. The following statements hold:
\begin{enumerate}
\item For all $n>0$ and all $a_0$, \dots, $a_{n-1}\in L$ such that
$\Lambda_n(a_0,\dots,a_{n-1})$ holds, there exists a unique
$\gf\colon\Co(n)\to L$ such that $\gf(\set{i})=a_i$, for all $i<n$.

\item For all $m$, $n>0$ and all $a_0$, $a_1$, \dots, $a_{m+n}\in L$
such that $\Lambda_{m,n}(a_0,a_1,\dots,a_{m+n})$ holds, there exists a
unique $\gf\colon L_{m,n}\to L$ such that $\gf(\set{i})=a_i$, for all
$i\neq m$, while $\gf(\set{m-1,m})=a_m$.
\end{enumerate}
\end{lemma}

\begin{proof}
Without loss of generality, $L$ is generated by $\setm{a_i}{0\leq i<n}$
in (i), by $\setm{a_i}{0\leq i\leq m+n}$ in (ii). In particular, by
Theorem~\ref{T:CLocFin}, $L$ is finite.
Since $L$ is a finite member of $\CC$,
we may assume, by Theorem~\ref{T:EmbThm}, that $L=\Co(T)$, for a finite
chain $(T,\utr)$. Let $u$ be the common value for
all $a_i\wedge a_j$ for $i\neq j$ in (i), for $i\neq j$ and
$\set{i,j}\neq\set{m-1,m}$ in (ii). The uniqueness statement about $\gf$
is, in both cases, obvious, and if there is a map $\gf$ as desired, then
it is given by the rule $\gf(X)=\bigvee_{i\in X}a_i$,
for all $X\in\Co(n)$ in (i), for all $X\in L_{m,n}$ in (ii), with
the convention that the empty join equals $u$. {}From
the assumption that the $a_i$-s satisfy (the statement involving joins
in) $\Lambda_n$ in (i) and $\Lambda_{m,n}$ in (ii) it follows easily that
$\gf$ is a \jh.

Now we prove that $\gf$ is a \mh.
Suppose first that $u$ is nonempty. The join of any two members of $L$ is
their union, whence $L$ is distributive. The statement that $\gf$ is a
\mh\ follows immediately in (i). In (ii), for all $X$, $Y\in L_{m,n}$,
we compute:
   \begin{align*}
   \gf(X)\wedge\gf(Y)&=\bigvee\famm{a_i\wedge a_j}{(i,j)\in X\times Y}\\
   &=\begin{cases}
   \gf(X\cap Y)\vee a_{m-1},&\text{if }
   (m-1,m)\in(X\times Y)\cup(Y\times X),\\
   \gf(X\cap Y),&\text{otherwise}.
   \end{cases}
   \end{align*}
But in the first case, $m-1$ belongs to $X\cap Y$, so we obtain again that
$\gf(X)\wedge\gf(Y)=\gf(X\cap Y)$.

Suppose now that $u=\es$. By Lemma~\ref{L:ChInt}, we may
assume without loss of generality that $a_i=[x_i,y_i)$, for elements
$x_i\utr y_i$ of $T$, for $i<n$ in (i) and $i\leq m+n$ in (ii), such that
   \[
   x_0\utr y_0\utr\cdots\utr x_{n-1}\utr y_{n-1}
   \]
in (i), while
   \[
   x_0\utr y_0\utr\cdots\utr x_{m-2}\utr y_{m-2}\utr x_m\utr y_m\utr
   \cdots\utr x_{m+n}\utr y_{m+n}
   \]
in (ii). Furthermore, from the assumption on the $a_i$-s it follows that
$x_{m-1}=x_m$ and $y_{m-1}\leq y_m$ in (ii), in particular,
   \[
   x_0\utr x_1\utr\cdots\utr x_{m+n}\text{ and }
   y_0\utr y_1\utr\cdots\utr y_{m+n}.
   \]
Hence, in both cases (i) and (ii), the value of $\gf(X)$ for $X$ in the
domain of $\gf$ can be computed by the rule $\gf(X)=[x_i,y_j)$ whenever
$X=[i,j]$, for $i\leq j$. It follows easily that $\gf$ is a \mh.
\end{proof}

Now we can prove the main result of the present section.

\begin{theorem}\label{T:SIProj}
Every finite subdirectly irreducible member of $\CC$ is projective in
$\CC$.
\end{theorem}

\begin{proof}
We first prove that $\Co(n)$ is projective in $\CC$, for all $n>0$.
Let $L\in\CC$, let $\pi\colon L\onto\Co(n)$ be a surjective lattice
homomorphism, we prove that there exists a lattice homomorphism
$\gf\colon\Co(n)\to L$ such that $\pi\circ\gf=\id_{\Co(n)}$. Since
$\Co(n)$ is finite, we may replace $L$ by a finitely generated
sublattice, which, by Theorem~\ref{T:CLocFin}, is finite. Since $L$ is
finite, the sublattice $\pi^{-1}\set{X}$ has a least element, that we
denote by $\beta(X)$, for any $X\in\Co(n)$. Put $a_i=\beta(\set{i})$,
for all $i<n$. Since $\beta$ is a \jh, the following statement holds:
   \begin{equation}\label{Eq:Ineqaijk}
   a_k\leq a_i\vee a_j,\text{ for }0\leq i<k<j<n.
   \end{equation}
Now we define inductively elements $b^l$ and $a_i^l$ of $L$, for $i<n$
and $l<\omega$, as follows:
   \begin{align}
   a_i^0&=a_i;\label{Eq:ai0}\\
   b^l&=\bigvee\famm{a_i^l\wedge a_j^l}{i\neq j};\label{Eq:bl}\\
   a_i^{l+1}&=a_i^l\vee b^l.\label{Eq:ail+1}
   \end{align}
Since $L$ is finite, there exists $l<\omega$ such that
$a_i^{l+1}=a_i^l$, for all $i<n$. {}From \eqref{Eq:Ineqaijk},
\eqref{Eq:ai0}, and \eqref{Eq:ail+1}, it is easy to prove, by induction
on $l$, the inequalities
   \begin{equation}\label{Eq:Ineqaijkl}
   a_k^l\leq a_i^l\vee a_j^l,\text{ for }0\leq i<k<j<n.
   \end{equation}
Furthermore, for $i\neq j$ and $k$ in $\fso{n-1}$,
$a_i^l\wedge a_j^l\leq b^l\leq a_k^{l+1}=a_k^l$.
Hence, the statement $\Lambda_n(a_0^l,\dots,a_{n-1}^l)$ holds, thus, by
Lemma~\ref{L:Lambdamn}(i), there exists a lattice homomorphism
$\gf\colon\Co(n)\to L$ such that $\gf(\set{i})=a_i^l$, for all $i<n$.
{}From \eqref{Eq:ai0}--\eqref{Eq:ail+1} it follows that
$\pi(a_i^l)=\set{i}$, for all $i<n$, whence $\pi\circ\gf=\id_{\Co(n)}$.

The proof that $L_{m,n}$ is projective, for $m$, $n>0$, is similar, by
using Lemma~\ref{L:Lambdamn}(ii). The definitions of the $a_i$-s and the
$a_i^l$-s are exactly the same as for the $\Co(n)$ case, while the join
in the definition of $b^l$ in \eqref{Eq:bl} has to be taken over the
pairs $(i,j)$ such that $i\neq j$ and $\set{i,j}\neq\set{m-1,m}$.

By Theorem~\ref{T:ClassSI}, there are no other finite subdirectly
irreducible members of $\CC$, thus the proof is complete.
\end{proof}

As a consequence of this, we obtain the following result, which shows
that $\CC$ is a quite peculiar variety, see the contrast with
Example~\ref{Ex:QnotV}.

\begin{theorem}\label{T:QVar2Var}
Every subquasivariety of $\CC$ is a variety.
\end{theorem}

\begin{proof}
Let $\QQ$ be a subquasivariety of $\CC$, we prove that
$\QQ$ is a variety. It suffices to prove that every homomorphic image
$L$ of a lattice $L'$ in $\QQ$ belongs to $\QQ$. Since $L$ belongs to the
locally finite variety $\CC$, it suffices to consider the case where $L$
is finite. By considering the subdirect decomposition of $L$, it
suffices then to consider the case where $L$ is subdirectly irreducible.
By Theorem~\ref{T:SIProj}, $L$ is projective within $\CC$, thus it
embeds into $L'$; whence $L$ belongs to $\QQ$.
\end{proof}

\section{An example}\label{S:ExpleNonPQ}

For a chain $Q$ and a subset $P$ of $Q$, endowed with the induced
ordering, the lattice $\Co(P)$ embeds into $\Co(Q)$, thus it belongs
to the variety generated by $\Co(Q)$. We shall now show, through an
example, that this simple observation cannot be extended to arbitrary
posets.

Let $P$ and $Q$ be the posets diagrammed on Figure~5. Obviously, $P$
is a subset of $Q$, endowed with the induced ordering.

\begin{figure}[htb]
\includegraphics{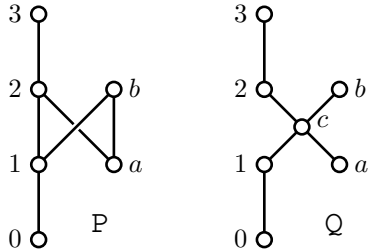}
\caption{The posets $P$ and $Q$}
\end{figure}

By induction on the natural number $n$, we define lattice terms
$x_1^{(n)}$, $x_2^{(n)}$, $s$, and~$t$, in
the variables $x_0$, $x_1$, $x_2$, $x_3$, $x_a$, $x_b$, putting
$x_1^{(0)}=x_1$, $x_2^{(0)}=x_2$, and
\begin{align*}
x_1^{(n+1)}&=
x_1^{(n)}\wedge\pI{x_0\vee x_2^{(n)}}\wedge(x_0\vee x_b);\\
x_2^{(n+1)}&=
x_2^{(n)}\wedge\pI{x_3\vee x_1^{(n)}}\wedge(x_3\vee x_a),
\end{align*}
for all $n<\omega$, then
$s=x_1\wedge\pII{x_0\vee\pI{(x_1\vee x_b)\wedge(x_2\vee x_a)}}$, and then
   \begin{multline*}
   t=(x_1\wedge x_b)
   \vee\pI{x_1\wedge(x_0\vee x_a)}\vee\pI{x_1\wedge(x_2\vee x_a)}\vee
   \pIII{x_1\wedge\pII{x_0\vee\pI{x_2\wedge(x_1\vee x_a)}}}\\
   \vee\pIII{x_1\wedge\pII{x_0\vee\pI{x_2\wedge(x_1\vee x_b)}}}
   \vee\pIII{x_1\wedge\pII{x_0\vee\pI{x_2\wedge(x_3\vee x_b)}}}.
   \end{multline*}
Finally, let $(*)$ be the following identity:
   \begin{equation}
   x_1^{(2)}\leq s\vee t.\tag{$*$}
   \end{equation}

\begin{lemma}\label{L:CoQhas*}
The lattice $\Co(Q)$ satisfies $(*)$.
\end{lemma}

\begin{proof}
Let $X_0$, $X_1$, $X_2$, $X_3$, $X_a$, $X_b$ be elements of $\Co(Q)$,
let $S$ and $T$ be obtained by evaluating $s$ and $t$ at those
parameters. We prove that $X_1^{(2)}$ is a subset of $S\cup T$. So,
let $x_1\in X_1^{(2)}$. If $x_1\in X_0\cup X_2$, then
$x_1\in(X_0\cap X_1)\cup(X_2\cap X_1)\subseteq T$; suppose now that
$x_1\notin X_0\cup X_2$. Since $x_1\in X_0\vee X_2^{(1)}$, there are
$x_0\in X_0$ and $x_2\in X_2^{(1)}$ such that either
$x_0\tr x_1\tr x_2$ or $x_2\tr x_1\tr x_0$.

Suppose that the first case occurs. If $x_2\in X_1\cup X_3$, then
   \[
   x_1\in\pII{X_1\cap\pI{X_0\vee(X_1\cap X_2)}}\cup
   \pII{X_1\cap\pI{X_0\vee(X_2\cap X_3)}}\subseteq T.
   \]
Suppose now that $x_2\notin X_1\cup X_3$.
Since $x_2\in X_1\vee X_3$, there are $x\in X_1$
and $x_3\in X_3$ such that either $x\tr x_2\tr x_3$ or
$x_3\tr x_2\tr x$. In the second case, from $x_1\tr x_2\tr x$ it follows
that $x_2\in X_1$, \contr. Thus $x_2\tr x_3$.

{}From $x_1\in X_1^{(2)}$ it follows that
$x_1\in X_0\vee X_b$. If $x_1$ belongs to $\dnw X_0$
(the lower subset of $Q$ generated by $X_0$), then
$x_1\in X_0\cap X_1\subseteq T$. If $x_1\in X_b$, then
$x_1\in X_1\cap X_b\subseteq T$. Suppose that
$x_1\notin\dnw X_0\cup X_b$. Since $x_1\in X_0\vee X_b$,
there exists $x_b\in X_b$ such that $x_1\tr x_b$.
Furthermore, from $x_2\in X_2^{(1)}$ it follows that
$x_2\in X_3\vee X_a$. If $x_2\in\dnw X_a$, then
$x_2\in X_2\cap(X_1\vee X_a)$, thus
$x_1\in X_1\cap\pII{X_0\vee\pI{X_2\cap(X_1\vee X_a)}}\subseteq T$.
Suppose now that $x_2\notin\dnw X_a$.
Since $x_2\notin X_3$ and $x_2\in X_3\vee X_a$, there exists $x_a\in X_a$
such that $x_a\tr x_2$.

If $x_a\utr x_1$, then $x_1\in X_1\cap(X_2\vee X_a)\subseteq T$.
If $x_1\utr x_a$, then $x_1\in X_1\cap(X_0\vee X_a)\subseteq T$.
Suppose now that $x_a\Vert x_1$ (where $\Vert$ denotes incomparability).
If $x_2\utr x_b$, then $x_2\in X_2\cap(X_1\vee X_b)$, thus
$x_1\in X_1\cap\pII{X_0\vee\pI{X_2\cap(X_1\vee X_b)}}\subseteq T$.
If $x_b\utr x_2$, then $x_2\in X_2\cap(X_3\vee X_b)$, thus
$x_1\in X_1\cap\pII{X_0\vee\pI{X_2\cap(X_3\vee X_b)}}\subseteq T$.
Suppose now that $x_2\Vert x_b$. Since $x_1\Vert x_a$, we have obtained
the inequalities
   \begin{equation}\label{Eq:ordxp}
   x_0\tr x_1\tr x_2\tr x_3,\qquad x_1\tr x_b,\qquad x_a\tr x_2,
   \qquad x_1\Vert x_a,\quad\text{ and }\quad x_2\Vert x_b.
   \end{equation}
This leaves the only possibility $x_p=p$, for all $p\in P$. In
particular,
   \[
   x_1=1\in\set{0}\vee\set{c}\subseteq
   X_0\vee\pI{(X_1\vee X_b)\cap(X_2\vee X_a)},
   \]
from which it follows that $x_1\in S$.

The other case to consider is $x_2\tr x_1\tr x_0$. Then, applying the
argument above to the dual of $\tr$, we obtain the dual of
\eqref{Eq:ordxp}, whence $x_k=3-k$, for all $k\in\set{0,1,2,3}$,
$x_a=b$, and $x_b=a$. In particular,
   \[
   x_1=2\in\set{3}\vee\set{c}\subseteq
   X_0\vee\pI{(X_1\vee X_b)\cap(X_2\vee X_a)},
   \]
from which it follows again that $x_1\in S$. In any case,
$x_1\in S\cup T$.
\end{proof}

\begin{lemma}\label{L:CoPnot*}
The lattice $\Co(P)$ does not satisfy $(*)$.
\end{lemma}

\begin{proof}
Put $x_p=\set{p}$, an element of $\Co(P)$, for any $p\in P$. Then the
left hand side of $(*)$, evaluated with those parameters, is
$x_1=x_1^{(2)}=\set{1}$, while the right hand side is empty.
Therefore, $\Co(P)$ does not satisfy $(*)$.
\end{proof}

Hence we have reached the desired conclusion.

\begin{proposition}\label{P:CoPnvCoQ}
The poset $P$ embeds into the finite poset $Q$, but
the lattice $\Co(P)$ does not belong to the variety generated
by~$\Co(Q)$.
\end{proposition}

\section{Open problems}\label{S:Pbs}

As in \cite{SeWe2}, we denote, for a class $\CK$ of posets, by
$\SUB(\CK)$ the lattice variety generated by
$\setm{\Co(P)}{P\in\CK}$. Say that a lattice variety $\VV$ is a
\emph{Stirlitz variety}, if it is of the form $\SUB(\CK)$ for some
class $\CK$ of posets.

It is clear that any join of Stirlitz varieties is a Stirlitz
variety, thus the set of all Stirlitz varieties, partially ordered by
inclusion, is a complete join-semilattice. In particular, it is a
lattice, however, we do not know whether the meet in this lattice is
the same as the meet for varieties.

\begin{problem}\label{Pb:IntSt}
Is the intersection of two Stirlitz varieties a Stirlitz variety?
\end{problem}

\begin{problem}\label{Pb:StCov}
Let $L$ be a lattice in $\SUB$. Does there exist a smallest Stirlitz
variety~$\VV$ such that $L\in\VV$?
\end{problem}

A related problem is the following.

\begin{problem}\label{Pb:FinMinSV}
For a finite lattice $L$ in $\SUB$, are there only finitely many
Stirlitz varieties $\VV$ which are minimal with the property that
$L\in\VV$?
\end{problem}

Analogies between our results with classical results of the
spatial theory of \emph{modular} lattices may fail. For example, the
main result of C. Herrmann, D. Pickering, and M. Roddy \cite{HPR}
states that every modular lattice embeds, \emph{within its
variety}, into an algebraic and spatial modular lattice. On the
other hand, every lattice $L$ in
$\SUB$ embeds into an algebraic and spatial lattice in
$\SUB$~--- namely, some $\Co(P)$, however, $\Co(P)$ may not belong to
$\Var(L)$, for example for $L=N_5$. This leads to the following
problem.

\begin{problem}\label{Pb:sameVar}
Does every lattice in $\SUB$ embed, \emph{within its variety}, into
some algebraic and spatial lattice?
\end{problem}

Of course, by Whitman's Theorem, every lattice $L$ embeds into a
partition lattice, which is both algebraic and spatial, but which does
not necessarily lie in the same variety as $L$. We do not even know
whether \emph{every} lattice embeds, within its variety, into an
algebraic and spatial lattice! While working on the present paper,
the authors met the following intriguing problem.

\begin{problem}\label{Pb:AlgDuAlg}
Can every lattice be embedded into some lattice that is both algebraic and
dually algebraic?
\end{problem}

\begin{noteadd}
The second author recently solved Problem~\ref{Pb:AlgDuAlg}.
\end{noteadd}

\begin{problem}\label{Pb:Sub(P)}
For a finite poset $P$, is the class of all sublattices of powers of
$\Co(P)$ a variety?
\end{problem}

The answer to Problem~\ref{Pb:Sub(P)} in the particular case
where $P$ is a \emph{chain} is, by the results of the
present paper, positive, see also Theorem~\ref{T:QVar2Var}. The results
of Section~\ref{S:ExpleNonPQ} also suggest a positive answer to
Problem~\ref{Pb:Sub(P)} in general.

\begin{example}\label{Ex:QnotV}

\begin{figure}[htb]
\includegraphics{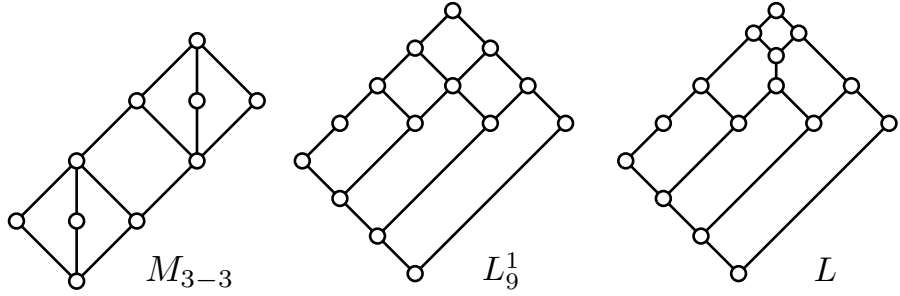}
\caption{The lattices $M_{3-3}$, $L_9^1$, and $L$}
\end{figure}

There are many finite lattices
$L$ for which the quasivariety $\QVar(L)$ generated by $L$ is
\emph{not} a variety, for example, the modular lattice $M_{3-3}$ of
Figure~6, see V.\,A. Gorbunov \cite[p.~257]{Gorb}. It is also
possible to find $L$ a bounded homomorphic image of a free lattice.
For example, the lattice $L_9^1$, see P. Jipsen and H. Rose
\cite{JiRo,JiRo2}, is bounded and subdirectly irreducible. It
also has a unique doubly reducible element; doubling this element
gives a finite, bounded lattice $L$. Furthermore, $L$ satisfies the
Whitman condition, thus it is projective, see \cite{FJN}. The lattices
$L_9^1$ and~$L$ are diagrammed on Figure~6.
Since $L_9^1$ is a quotient of $L$, it belongs
to $\Var(L)$. If $L_9^1$ belonged to
$\QVar(L)$, then, since it is subdirectly
irreducible, it would embed into~$L$, which is easily seen not to be
the case. Therefore, $\QVar(L)\neq\Var(L)$. Compare this with
Theorem~\ref{T:QVar2Var}.

\end{example}

\begin{problem}\label{Pb:ConLattSUB}
What are the congruence lattices of lattices in $\CC$?
\end{problem}

Our next problems are related to the variety $\Csub(\bD)$ studied by
V. Slav\'\i k in \cite{Slav}. This variety contains the variety $\CC$
studied in the present paper, see the proof of Theorem~\ref{T:CLocFin}. In
\cite{Slav}, some properties of the finite subdirectly irreducible
members of $\Csub(\bD)$ are given, for example, every proper dual ideal
is a distributive lattice.

\begin{problem}\label{Pb:SICsubD}
Describe the lattice of subvarieties and classify the finite subdirectly
irreducible members of $\Csub(\bD)$.
\end{problem}

In V. Slav\'\i k \cite{Slav2}, it is proved that $\Csub(\bD)$ has
uncountably many subvarieties, but this does not seem to rule out a
reasonable classification of finitely generated subvarieties.

Furthermore, it is proved in \cite{Slav} that $\Csub(\bD)\cap\bM=\bMo$,
where $\bM$ (resp.,~$\bMo$) denotes the variety of all modular lattices
(resp., the variety generated by the infinite countable lattice
$M_\omega$ of length two). It is well-known that $\bMo$ is finitely
based, see B.~J\'onsson~\cite{Jons68} or \cite[Theorem~3.32]{JiRo}. This
suggests the following problems.

\begin{problem}\label{Pb:CsubFB}
Is the variety $\Csub(\bD)$ finitely based? More generally, if $\VV$ is
a self-dual, finitely based variety of lattices, is $\Csub(\VV)$
finitely based?
\end{problem}

\begin{problem}\label{Pb:IntersSD+}
Describe $\Csub(\bD)\cap\mathbf{SD}_\vee$, where
$\mathbf{SD}_\vee$ denotes the quasivariety of all \jsd\ lattices. In
particular, is $\Csub(\bD)\cap\mathbf{SD}_\vee$ a finitely based
quasivariety?
\end{problem}

\section*{Acknowledgment}

This work was completed during the first author's visit at the
University of Caen in March and April 2002, supported by a Young
Scientist INTAS fellowship program. The hospitality of the SDAD team can
never be forgotten.

This work was started during the two authors' visit at the Charles
University, from October to December 2001. Special thanks
are due to Ji\v{r}\'\i\ T\r{u}ma and V\'aclav Slav\'\i k.

\end{document}